\documentclass[reqno, 12pt]{amsart}
\usepackage{amssymb,amsmath,amsfonts}
\usepackage[margin=1in]{geometry}
\usepackage{dsfont}
\usepackage{color}
\usepackage{graphicx,graphics}
\usepackage{latexsym}
\usepackage[dvips]{epsfig}
\usepackage{mathtools,mathrsfs}
\usepackage{stmaryrd,cite}
\usepackage{cases}
\usepackage{enumerate}
\usepackage[usenames,dvipsnames,svgnames]{xcolor}
\usepackage{todonotes}
\usepackage{marginnote}
\usepackage{listings}
\usepackage{cleveref}
\allowdisplaybreaks
\definecolor{refkey}{rgb}{1,0,0.5}
\definecolor{labelkey}{rgb}{0,0.4,1}
\numberwithin{equation}{section}
\newtheorem{thm}{Theorem}[section]

\newtheorem{lem}[thm]{Lemma}
\newtheorem{prop}[thm]{Proposition}
\newtheorem{rmk}[thm]{Remark}

\newcommand{\ea}{\epsilon}
\newcommand{\vea}{\varepsilon}

\newcommand{\ka}{\kappa}
\newcommand{\al}{\alpha}

\newcommand{\da}{\delta}

\newcommand{\la}{\lambda}
\newcommand{\na}{\nabla}

\newcommand{\ta}{\theta}

\newcommand{\sa}{\sigma}

\newcommand{\ga}{\gamma}
\newcommand{\ba}{\beta}

\newcommand{\oa}{\omega}
\newcommand{\Oa}{\Omega}

\newcommand{\iy}{\infty}
\newcommand{\pl}{\partial}

\newcommand{\lt}{\left}
\newcommand{\rt}{\right}
\newcommand{\be}{\begin{equation}}
\newcommand{\ee}{\end{equation}}
\newcommand{\bee}{\begin{equation*}}
\newcommand{\eee}{\end{equation*}}

\newcommand{\ef}{\eqref}
\newcommand{\f}{\frac}
\newcommand{\les}{\lesssim}

\begin{document}
\title[] {Global Existence and Asymptotic Equivalence to Barenblatt-type Solutions for the Physical Vacuum Free Boundary Problem  of Damped Compressible Euler Equations in M-D}
\author{Huihui Zeng}
\maketitle

\begin{abstract}
For the physical vacuum free boundary problem of the damped compressible Euler equations in both 2D and 3D, we prove the global existence of smooth solutions and justify their time-asymptotic equivalence to the corresponding Barenblatt self-similar solutions derived from the porous media  equation under Darcy's law approximation, provided the initial data are small perturbations of the Barenblatt solutions.
Building on the 3D almost global existence result in [Zeng, Arch. Ration. Mech. Anal. 239, 553--597 (2021)], our key contribution lies in improving the decay rate of the time derivative of the perturbation from $-1$ (as previously established) to $-1-\varepsilon$ for a fixed constant $\varepsilon > 0$. This critical enhancement ensures time integrability and hence global existence. Together with the previous 1D result in [Luo--Zeng, Comm. Pure Appl. Math. 69, 1354--1396 (2016)],
the results obtained in this paper provide a complete answer to the question raised in [Liu, T.-P.: Jpn. J. Appl. Math. 13, 25--32 (1996)]. Moreover, we also consider the problem with time-dependent damping of the form $(1+t)^{-\lambda}$ for $0 < \lambda < 1$. Notably, our framework unifies the treatment of both time-dependent ($0 < \lambda < 1$) and time-independent ($\lambda = 0$) damping cases across dimensions.
We further quantify the decay rates of the density and velocity, as well as the expansion rate of the physical vacuum boundary.
\end{abstract}

\section{Introduction}
Consider the following vacuum free boundary problem for isentropic compressible Euler equations with  damping:
\begin{subequations}\label{2.1} \begin{align}
& \pl_t \rho   + {\rm div}(\rho   u ) = 0 &  {\rm in}& \ \ \Omega(t), \label{2.1a}\\
 &  \pl_t  (\rho   u )   + {\rm div}(\rho   u \otimes   u )+\nabla  p(\rho) = - (1+t)^{-\la}\rho {  u}  & {\rm in}& \ \ \Omega(t),\label{2.1b}\\
 &\rho>0 &{\rm in }  & \ \ \Omega(t),\label{2.1c}\\
 & \rho=0    &    {\rm on}& \  \ \Gamma(t)=\pl \Omega(t), \label{2.1d}\\
 &    \mathcal{V}(\Gamma(t))={  u}\cdot \mathcal{N}, & &\label{2.1e}\\
&(\rho,{ u})=(\rho_0, { u}_0) & {\rm on} & \ \   \Omega(0). \label{2.1f}
 \end{align} \end{subequations}
Here $( t,x)\in [0,\iy)\times\mathbb{R}^n $,  $\rho $, ${ u} $, and $p$ represent, respectively, the time and space variable, density, velocity and  pressure; $\Omega(t)\subset \mathbb{R}^n$, $\Gamma(t)$, $\mathcal{V}(\Gamma(t))$ and $ \mathcal{N}$ denote, respectively, the changing volume occupied by the gas at time $t$, moving vacuum boundary, normal velocity of $\Gamma(t)$, and exterior unit normal vector to $\Gamma(t)$;
$\la \in [0,1)$ is a constant describing the time decay exponent of the damping coefficient, with the special case of $\lambda=0$ corresponding to the classical time-independent damping.
We hypothesize that the pressure satisfies the $\gamma$ law:
$$ p(\rho)=\rho^{\gamma}, \  \ {\rm where} \ \  \gamma>1 {~\rm is~the~adiabatic~exponent}; $$
and the initial density satisfies
\begin{subequations}\label{initial density}\begin{align}
& \rho_0>0 \ \ {\rm in} \ \ \Omega(0), \ \  \rho_0=0 \ \ {\rm on} \ \ \Gamma(0),   \ \ \int_{\Omega(0)} \rho_0(x) dx =M,
\\
&  -\infty<\nabla_\mathcal{N}\lt(c^2(\rho_0)\rt)<0  \  \ {\rm on} \ \  \Gamma(0), \label{4.21-1}
\end{align}\end{subequations}
where $M\in (0, \iy)$ is the initial total mass, and $c(\rho)=\sqrt{ p'(\rho)}$ is the sound speed.
Indeed, condition \ef{4.21-1} defines  a physical vacuum boundary or singularity (cf. \cite{7,10',16',23,24,25}).

The compressible Euler equations of isentropic flows with damping, \ef{2.1a}-\ef{2.1b}, are closely related to the
porous media type equation (cf. \cite{HL,HMP,HPW, 23,LZ,HZ}):
\begin{equation}\label{pm}
\pl_t \rho  = (1+t)^{\la}\Delta p(\rho),
\end{equation}
when the momentum equation \eqref{2.1b} is simplified to the Darcy-type law:
\begin{equation}\label{darcy}
\nabla  p(\rho)=- (1+t)^{-\la} \rho { u}.
\end{equation}
A fundamental class of solutions to equation \eqref{pm} with finite total mass
$M$ is provided by the self-similar solutions:
\be\label{Feb25-1}
\bar \rho(t,x)=\nu^{-n}(t)
\left(\bar{A}- \bar{B} \nu^{-2}(t)|x|^2\right)^{1/({\gamma-1})}.
\ee
defined in the expanding region
$
\bar \Omega(t)= \{x\in\mathbb{R}^n:\  |x| < \nu(t) \sqrt{\bar{A}/ \bar{B}} \}
$,
where
\be\label{5.26-1}
\nu(t)=(1+t)^{\ka}, \ \  \ka=\frac{1+\lambda}{n\ga-n+2},  \ \   \bar{B}=\frac{\ga-1}{2\ga}\ka,
\ee
and  $\bar{A}$ is a  positive constant determined by $n$, $\la$, $\ga$ and $M$. In the special case when $\la=0$, solution \eqref{Feb25-1} reduces to the celebrated Barenblatt self-similar solution for the porous media equation (cf. \cite{ba}).
The Barenblatt-type solution in  \ef{Feb25-1} satisfies $\bar\rho> 0$ in $\bar \Omega(t)$, $ \bar\rho=0$ on $ \pl \bar \Omega(t)$, and $\int_{\bar \Omega(t)}  \bar\rho(t,x )dx  =M$
for all  $t\ge 0$. We assume that the initial total mass of problem \ef{2.1} coincides with that of \ef{pm}:
$$
\int_{\Omega(0)} \rho_0(x) dx =M=\int_{\bar \Omega(t)}  \bar\rho(t,x )dx.
$$
Clearly, $\pl \bar\Oa(t)$ represents a physical vacuum boundary, since $\nabla_\mathcal{N}\lt(c^2(\bar\rho)\rt)=-2\ga \sqrt{\bar A \bar B}(1+t)^{\ka-1-\la}$ on $\pl \bar\Oa(t)$.
The corresponding  Barenblatt-type  velocity field $\bar { u}$ is defined by
\be\label{4.24-4}
\bar { u}(t,x)=-   (1+t)^{\la}  \f{ \nabla p(\bar \rho) }{\bar\rho}
=\ka \f{  x }{1+t}
 \ \ {\rm in} \ \  \bar \Omega (t),
 \ee
and thus the system
\ef{pm}-\ef{darcy} is solved by
$(\bar\rho, \bar { u})$  in  $\bar \Omega(t)$.

It is important to understand the physical vacuum of the  free boundary problem \eqref{2.1} and to study the time-asymptotic behavior of its solutions in relation to the Barenblatt-type self-similar solution. This line of research was initiated  in \cite{23}. Specifically, for problem \eqref{2.1} with time-independent damping ($\lambda=0$ in \eqref{2.1}),  a class of particular solutions was constructed in  \cite{23}, using the following ansatz:
$$
\Omega(t) = B_{R(t)}(0), \quad c^2(x, t) = e(t) - b(t) r^2, \quad u(x, t) = \left( \frac{x}{r} \right) u(r, t),
$$
where $r = |x|$, $R(t) = \sqrt{e(t)/b(t)}$, and $u(r, t) = a(t) r$.
A system of ordinary differential equations for $(e, b, a)(t)$ was derived in \cite{23}, with $e(t), \ b(t) > 0$ for $t \ge 0$. This family of particular solutions was shown to be time-asymptotically equivalent to the Barenblatt self-similar solution with the same total mass.
Since the construction of these particular solutions for $\lambda=0$ in \cite{23}, a natural and important open question has been whether, for general initial data, a global-in-time existence theory exists for problem \eqref{2.1} that captures the physical vacuum  behavior \eqref{4.21-1}, and whether the time-asymptotic equivalence still holds between solutions of \eqref{2.1} and the corresponding Barenblatt self-similar solution with the same total mass.
This question has been answered affirmatively in the 1D case (cf. \cite{LZ}) and in the 3D spherically symmetric case (cf. \cite{HZ}). For general 3D perturbations without symmetry assumptions, the author of the present paper constructed an {\it almost} global-in-time solution in \cite{HZeng}. Note that the {\it almost} global result in \cite{HZeng} is still finite in nature; thus, extending it from finite to infinite time is critical both in the theory of nonlinear PDEs and in applications.
The first aim of this paper is to extend the almost global result in \cite{HZeng} to a global  one for the time-independent damping case ($\lambda=0$) in both 2D and 3D, and to justify the time-asymptotic equivalence to the Barenblatt self-similar solutions. This provides a complete answer to the open question raised in \cite{23}, building on a series of efforts in \cite{LZ, HZ, HZeng}.
Moreover, we also consider the problem with time-dependent damping, where the damping coefficient is $(1+t)^{-\lambda}$ with $0 < \lambda < 1$. Notably, our framework unifies the treatment of both time-dependent ($0 < \lambda < 1$) and time-independent ($\lambda = 0$) damping cases across dimensions.
We further quantify the decay rates of the density and velocity, as well as the expansion rate of the physical vacuum boundary.

To contextualize our approach, we first review the relevant literature on the long-time existence theory for physical vacuum free boundary problems, in particular,
\cite{HZeng}, in which the almost global existence result was obtained for problem \eqref{2.1}-\eqref{initial density} in the case of  time-independent damping in 3D.
Using the formulation in \cite{HZeng},  we define the background Lagrangian flow as $\tilde x(t,y )=\lt( \nu(t) +   h(t) \rt) y $, where the higher-order correction $h(t)$ is  determined by an initial value problem of ODEs such that $\nu+h$ maintains the same asymptotic behavior as $\nu$.
The role of $h(t)$ is to compensate for the error arising from the fact that the Barenblatt-type solution solves the porous media type equation, not the original compressible equations with damping.
It is then natural to decompose $x(t, y)$, the Lagrangian flow of the velocity field $u$,  as
$$x(t,y)=\tilde x(t,y ) +(\nu+h)(t) \omega(t,y)=(\nu+h)(t) \lt(y+\omega(t,y)\rt),$$
reducing  the problem to studying the equations for the perturbation $\omega(t,y)$.
The higher-order energy norms defined in \cite{HZeng} take the form
\be\label{a} \mathscr{E}(t)=\mathscr{E}_1(t)
+ \nu^{2\ka^{-1}}(t)\mathscr{E}_2(t)=
\mathscr{E}_1(t)
+ (1+t)^{2}\mathscr{E}_2(t),\ee
where $\mathscr{E}_1(t)$ and  $\mathscr{E}_2(t)$ represent the squared weighted $L^2$-norms of the spatial derivatives of $\omega$ and the space-time mixed derivatives  of $\pl_t \omega$, respectively.
This  energy norm structure was
employed  in \cite{LZ}  to establish the global existence for  problem \ef{2.1}-\ef{initial density} with time-independent damping  in 1D,  and  in \cite{HZ} for  the  spherically symmetric 3D case, where the key  global existence estimate takes the form
\be\label{7.15}
 \mathscr{E}(t) \les  \nu^{-2}(t)  \mathscr{E}(0)=(1+t)^{-2\ka} \mathscr{E}(0).
 \ee
Here and thereafter, we use the notation $a\les b$ to denote $a\le C b$ for some positive constant $C$ independent of $t$.
However,  obtaining \ef{7.15}  becomes  practically impossible for the general  3D case without symmetry assumptions, due to the presence of the curl of $\oa$. Indeed, as shown in \cite{HaJa1,ShSi},   the curl of $\oa$ evolves as
\be\label{7.14}
{\rm curl} \oa(t,y)= {\rm curl} \oa(0,y) + \int_0^t  {\rm curl} \pl_\tau \oa(\tau, y) d\tau.
\ee
Consequently, the sharpest possible estimate for $\mathscr{E}_1(t) $ is limited to
\be\label{7.17}
\mathscr{E}_1(t) \les  \mathscr{E}(0),
\ee
since the estimate of $\mathscr{E}_1(t)$ involves ${\rm curl} \oa$ in the general asymmetric case.
The expression \eqref{7.14} is essential  for curl estimation because while no direct equation exists for  ${\rm curl} \oa$,  the quantity  ${\rm curl} \pl_t \oa$  satisfies an equation derived from \ef{2.1b} that is more tractable for analysis.
Since previous studies of the global existence for physical vacuum free boundary problems (cf. \cite{LZ,HaJa1,ShSi, HZ}) estimate $\mathscr{E}_1(t)$ and $\mathscr{E}_2(t)$  simultaneously, and
the temporal decay of $\pl_t \omega$  and its derivatives is naturally encoded in the definition of $\mathscr{E}(t)$,   the estimate achieved  in \cite{HZeng} yields
\be\label{7.16}
\mathscr{E}(t) \les (\ln(1+t)+1) \mathscr{E}(0).
\ee
The logarithmic term in \eqref{7.16} arises from  \eqref{7.14} because ${\rm curl} \partial_t \omega$ decays asymptotically like $(1+t)^{-1}$.
This estimate  explains why the result in \cite{HZeng}was  limited to the almost global existence rather than full global existence.

We now outline the key differences between our framework and \cite{HZeng}, emphasizing the methodological advances that yield the global existence.
As established in \eqref{7.14}, the global existence hinges critically on the integrability of $\operatorname{curl} \partial_t \omega$ over all time.
To ensure this integrability, our approach diverges fundamentally from prior works (cf. \cite{LZ,HaJa1,ShSi,HZ,HZeng}), which treated $\mathscr{E}_1(t)$ and $\mathscr{E}_2(t)$ as a unified quantity in their estimates. Instead, we treat $\mathscr{E}_1(t)$ and $\mathscr{E}_2(t)$ as distinct but coupled quantities, exploiting their individual behaviors for sharper bounds.
While \eqref{7.17} suggests that temporal decay for $\mathscr{E}_1(t)$ is unlikely, such decay remains attainable for $\mathscr{E}_2(t)$ and by extension, for $\nu^{2\kappa^{-1}}(t)\mathscr{E}_2(t)$, as evidenced by \eqref{a} and \eqref{7.15}.  In spherical symmetry (cf. \cite{HZ}), one indeed obtains
$$\nu^{2\kappa^{-1}}(t)\mathscr{E}_2(t) \les \nu^{-2}(t)  \mathscr{E}(0).$$
However, the inherent complexity of asymmetric flows precludes this ideal decay rate in the general 3D case. Our analysis reveals that the sharpest achievable bound takes the form
$$\nu^{2\kappa^{-1}}(t)\mathscr{E}_2(t) \les \nu^{-2+\varepsilon}(t)  \mathscr{E}(0)$$
for an arbitrary positive constant $\varepsilon$.
To streamline the presentation, we introduce the following higher-order energy norm for this work:
\begin{equation}\label{7.20}
\mathscr{E}(t)=\mathscr{E}_1(t)
+ \nu^{2\kappa^{-1}+1}(t)\mathscr{E}_2(t)=
\mathscr{E}_1(t)
+ (1+t)^{2+\kappa}\mathscr{E}_2(t),
\end{equation}
which differs from \ef{a}.
Crucially, we prove the uniform bound:
\begin{equation}\label{7.18}
 \mathscr{E}(t)\les   \mathscr{E}(0),
\end{equation}
which establishes the global existence of smooth solutions to \eqref{2.1}-\eqref{initial density}.

To establish \eqref{7.18}, we outline the key steps of our analysis, which represent a significant departure from previous approaches in \cite{LZ,HaJa1,ShSi,HZ,HZeng}.
While prior works focused primarily on the equation for $\omega$, our analysis instead centers on the time-differentiated equation governing $\partial_t \omega$.
The system, as discussed in \cite{LZ,HZ,HZeng}, exhibits dual hyperbolic-parabolic behavior: the hyperbolic nature dominates when the damping term is treated as a lower-order perturbation, while degenerate parabolic behavior emerges when the acceleration term decays rapidly, consistent with Darcy's law at large times.
For the spherically symmetric case \cite{HZ}, the proof of \eqref{7.15} employs weighted hyperbolic-parabolic estimates for both the $\omega$-  and $\partial_t \omega$-equations. However, in asymmetric flows, the situation differs between the equations: parabolic estimates fail for the $\omega$-equation (as evidenced in \eqref{7.17}), whereas both hyperbolic and parabolic estimates remain valid for the $\partial_t \omega$-equation.
Our proof strategy proceeds as follows:
first estimate $\mathscr{E}_2(t)$ by leveraging the parabolic structure of the $\partial_t \omega$  equation;
then derive bounds for $\mathscr{E}_1(t)$ based on the $\mathscr{E}_2(t)$ estimates.
This approach requires a careful decoupling of $\mathscr{E}_1(t)$ and $\mathscr{E}_2(t)$ in the analysis, particularly in handling nonlinear terms that depend on both quantities.
These refined techniques extend beyond the methods in \cite{HZ,HZeng}, providing the framework to overcome the challenges of asymmetric flows and rigorously establish the global existence.

We now turn to the discussion of the  time-dependent damping mechanisms. These present three key analytical challenges that do not arise in the time-independent scenario: (i) parametric complexity: the time decay exponent $\lambda$ of the damping coefficient introduces additional complexity in the analysis; (ii) asymptotic weakening: the damping effect diminishes as $t\to\infty$; (iii) dual behavior of correction term: the term $h(t)$ exhibits algebraic decay in the time-independent case but switches to algebraic growth in the time-dependent case when  $\lambda>1-\kappa$.

For the 3D physical vacuum free boundary problem of the compressible Euler equations, the available global-in-time results \cite{HaJa1,ShSi} are confined to expanding solutions where the vacuum boundary expands linearly at a rate of $O(1+t)$, arising from initial data that are small perturbations of affine motions \cite{sideris1,sideris2}. (See also \cite{CHJ,HaJa2,PHJ} for related results on global expanding solutions with linear expansion rates.) The stabilizing influence of fluid expansion is also relevant in other contexts, such as in general relativistic cosmological models \cite{Rodnianski,Oliynyk}.
A distinctive feature of our problem \eqref{2.1}, which sets it apart from the aforementioned global results, is the \emph{sub-linear} expansion of its vacuum boundary. To illustrate this, consider the time-independent damping case ($\lambda=0$). Here, the Barenblatt solutions---which serve as the long-time background profiles---exhibit a vacuum boundary that expands at a sub-linear rate, specifically $O((1+t)^{1/(3\gamma-1)})$ in 3D. For $\gamma > 1$, this rate is slower than $O((1+t)^{1/2})$. This slow expansion introduces significant difficulties in proving global existence for solutions to \eqref{2.1}, primarily due to the gradual decay of various quantities.
In this work, we overcome these challenges and prove the global existence of smooth solutions to problem \eqref{2.1} for initial data that are small perturbations of Barenblatt-type self-similar solutions.
For the M-D physical vacuum free boundary problem of inviscid compressible fluids with sub-linear boundary expansion, this result establishes, for the first time and without any symmetry assumptions, the global existence of solutions.

The study of vacuum states in compressible fluids has a long history, dating back to \cite{LiuSmoller}, which demonstrated that shock waves vanish at the vacuum. Early well-posedness results were established for sound speeds smoother than $C^{1/2}$-H\"older continuous \cite{chemin1,chemin2,24,25,MUK,Makino,38,39}. The physical vacuum, where the sound speed $c(\rho)$ is exactly $C^{1/2}$-H\"older continuous, poses greater challenges since the standard theory of symmetric hyperbolic systems \cite{Friedrichs,Kato,17} fails. More recent advances have established local well-posedness for the compressible Euler equations with physical vacuum \cite{16,10,7,10',16', IfTa, LiuL} (see also \cite{zhenlei,zhenlei1,LXZ,serre}).
Physical vacuum arises naturally in other important physical contexts such as gaseous stars \cite{6',cox,HaJa2,LXZ}. A key motivation is to understand the long-time stability of physically significant profiles with scaling invariance such as self-similar solutions, affine motions, and equilibria¡ªwhich require high regularity near the vacuum boundary. Extending local solutions to global ones is particularly challenging due to the strong degeneracy. For the 1D Cauchy problem, $L^p$-convergence of $L^\infty$ weak solutions to Barenblatt profiles was shown in \cite{HMP,HPW} via entropy estimates; however, such methods do not track the vacuum interfaces and are difficult to apply in M-D settings without symmetry assumptions.
In the context of compressible Euler equations with time-dependent damping, the studies in \cite{HY, HWY} established global existence and blow-up results for smooth solutions in M-D for the Cauchy problem away from vacuum states. Related results on physical vacuum free boundary problems in 1D and the 3D spherically symmetric case were obtained in \cite{Panxh1,Panxh2}.

The remainder of the paper is structured as follows. Section~\ref{sec2} reformulates the vacuum free boundary problem in Lagrangian coordinates, introduces the perturbation framework around the Barenblatt-type self-similar solutions, and states the main global existence and stability results in Theorems~\ref{mainthm} and~\ref{orig}.
Section~\ref{sec3} proves Theorem~\ref{mainthm} by establishing global existence through a priori estimates. It begins with preliminary estimates, followed by curl estimates (Lemma~\ref{prop-curl}) and energy estimates (Lemma~\ref{3.31-1}). The core of the argument involves a decoupled analysis of the energy functionals \(\mathscr{E}_1(t)\) and \(\mathscr{E}_2(t)\), along with a refined treatment of the time-differentiated equation for \(\partial_t \omega\).
Section~\ref{sec4} presents the proof of Theorem~\ref{orig}, quantifying the convergence rates of the density and velocity to their Barenblatt-type profiles, as well as the long-time behavior of the vacuum boundary. Finally, the Appendix compiles supplementary technical tools.

\section{Reformulation of the problem and main results}\label{sec2}

We set the initial domain $\bar \Omega(0)$ as the reference domain $\Oa$ by $\Oa= \{x\in\mathbb{R}^n:\  |x| < \sqrt{\bar{A}/ \bar{B}} \}$,  define $x$ as the Lagrangian flow of the velocity $u$ by
$$
\pl_t x(t,y)= u(t, x(t,y)) \  {\rm for} \  t>0,  \  {\rm and} \   x(0,y)=x_0(y) \ {\rm for} \   y\in \Omega,
$$
and  set the Lagrangian density,  the inverse of the Jacobian matrix, and the Jacobian determinant by
$$
\varrho(t,y)=\rho(t, x(t,y)),  \ \
\mathscr{A}(t,y)=\lt(\frac{\pl x}{\pl y}\rt)^{-1}, \ \
\mathscr{J}(t,y)={\rm det} \lt(\frac{\pl x}{\pl y}\rt)  .
$$
Then,  system \ef{2.1} can be written in Lagrangian coordinates as
\begin{subequations}\label{2.1n}\begin{align}
& \pl_t  \varrho  +   \varrho \mathscr{A}_i^k \pl_t\pl_k x^i = 0 &  {\rm in}& \ \ (0,\iy)\times\Oa, \label{2.1na}\\
 &   \varrho \pl_{tt} x_i   + \mathscr{A}_i^k \pl_k ( { \varrho^\ga} )  = -(1+t)^{-\la}  \varrho \pl_t x_i  & {\rm in}& \ \ (0,\iy)\times\Oa,\label{2.1nb}\\
 & \varrho>0 &{\rm in }  & \ \ (0,\iy)\times\Oa
 ,\label{2.1nc}\\
 &  \varrho=0    &    {\rm on}& \  (0,\iy)\times\pl\Oa, \label{2.1nd}\\
&( \varrho,   x, \pl_t x)=(\rho_0(x_0), x_0, u_0(x_0)) & {\rm on} & \ \  \{t=0\}\times\Omega, \label{2.1nf}
 \end{align} \end{subequations}
where $x^i=x_i$ and $\pl_k=\frac{\pl}{\pl y_k}$. It follows from  \ef{2.1na} and
$\pl_t \mathscr{J} = \mathscr{J} \mathscr{A}_i^k \pl_t\pl_k x^i   $ that
$$ \varrho(t,y)\mathscr{J}(t,y)= \varrho (0,y)\mathscr{J}(0,y)=\rho_0\lt(x_0(y)\rt) {\rm det} \lt(\frac{\pl x_0(y)}{\pl y}\rt) .$$
We choose $x_0(y)$ such that
$\rho_0\lt(x_0(y)\rt) {\rm det} (\frac{\pl x_0(y)}{\pl y})=   \bar\rho(0,y)$. The existence of such an $x_0$ follows from the Dacorogna-Moser theorem (cf. \cite{DM}) and \ef{initial density}. So the Lagrangian density can be expressed as
\be\label{4.24-1}
\varrho=\bar\rho_0 \mathscr{J}^{-1}, \ \ {\rm where} \ \
\bar\rho_0(y) = \bar\rho(0,y)= \lt( \bar{ A }-  \bar{B} |{y}|^2 \rt)^{{1}/({\ga-1})},
\ee
and problem \ef{2.1n} reduces to
\begin{subequations}\label{system}\begin{align}
&\bar\rho_0 \pl_{tt}  x_i    + \mathscr{J}  \mathscr{A}_i^k\pl_k  \lt(\bar\rho_0^\ga \mathscr{J}^{-\ga}\rt) = -  (1+t)^{-\la} \bar\rho_0  \pl_t x_i  &  {\rm in}& \ \  (0,\iy)\times\Omega, \label{171018} \\
&(   x, \pl_t x)=(x_0, u_0(x_0)) & {\rm on} & \ \   \{t=0\} \times \Omega.
\end{align}\end{subequations}

We define $\bar x$ as the Lagrangian flow of the   velocity $\bar u$ by
$$\pl_t \bar x(t,y)= \bar u(t, \bar x(t,y))  \ {\rm for} \  t>0, \  {\rm and} \    \bar x (0,y) =y  \  {\rm for} \   y\in \Omega,$$
thus,
\be\label{4.24-2}
\bar x(t,y)=\nu(t) y \ \ {\rm for} \ \  y\in  \Omega , \ \  {\rm where} \ \
\nu(t)=(1+t)^{\ka }.
\ee
Since $\bar x(t, y)$ does not solve equation \ef{171018}, we introduce a correction $h(t)$, as in \cite{LZ, HZ,HZeng},  which is the  solution to the following initial value problem of ordinary differential equations:
\begin{subequations}\label{pomt}\begin{align}
& h_{tt} +  { (1+t)^{-\la} } h_t -\ka (\nu+h)^{n-n\ga-1 }=-  \nu_{ tt}  - { (1+t)^{-\la} }\nu_{t} , \ \ t >0,  \\
& h(t=0)=h_t(t=0)=0.
 \end{align}\end{subequations}
It is worth noting that $\nu+h$ exhibits qualitative behavior similar to $\nu$, see \ef{5.26} for details.
The ansatz is then given by
$$
\tilde x(t,y )=\bar x (t,y )  +   h(t) y =\ta(t) y,  \ {\rm where}  \ \ta(t)=\nu(t)+h(t),
$$
which satisfies
$$
\bar\rho_0  \pl_{tt} \tilde x_i   +  \tilde{\mathscr{J}}\tilde{\mathscr{A}}_i^k  \pl_k\lt(\bar\rho_0^\ga \tilde{\mathscr{J}}^{-\ga}\rt)  = -  (1+t)^{-\la}\bar\rho_0  \pl_t \tilde x_i \ {\rm in} \  (0,\iy) \times \Omega,
$$
where $\tilde{\mathscr{J}}= {\rm det} (\frac{\pl \tilde x}{\pl y})=\ta^n$, $\tilde{\mathscr{A}}=(\frac{\pl \tilde x}{\pl y})^{-1} = \ta^{-1}\mathbb{I}$, and
$\mathbb{I}$ is the identity matrix.

We write the problem in the perturbation form around the ansatz $\tilde x(t,y )$ defined above. We set
\begin{align*}
& \eta(t,y)=\ta^{-1}(t)x(t,y), \  \  \oa(t,y)  = \eta(t,y)-y , \\
& A(t,y)  =\lt(\frac{\pl \eta}{\pl y}\rt)^{-1}=\lt(\mathbb{I}+\frac{\pl \oa}{\pl y}\rt)^{-1}, \
 J(t,y)  = {\rm det} \lt(\frac{\pl \eta }{\pl y}\rt) ={\rm det} \lt(\mathbb{I}+\frac{\pl \oa }{\pl y}\rt),
\end{align*}
then problem \ef{system}, hence problems \ef{2.1n} and \ef{2.1}, can be written as
\begin{subequations}\label{newsystem}\begin{align}
&\ta^{n\ga-n+2} \bar\rho_0   \pl_t^2 \oa_i + \ta^{n\ga-n+1} \lt ( (1+t)^{-\la}\ta+ 2\ta_t \rt) \bar\rho_0 \pl_t  \oa_i \notag\\
&\quad + \ka \bar\rho_0 \oa_i    +   \pl_k \lt(  \bar\rho_0^\ga  (A_i^k {J}^{1-\ga} -  \da_i^k ) \rt)  = 0 \ \ &  {\rm in} & \ \ (0,\iy)\times\Omega, \label{3-1-3}\\
& (  \oa , \pl_t \oa)=\lt ( x_0-y, u_0(x_0)- \ka x_0\rt)  \ \  & {\rm on}  \ \  & \{t=0\}\times \Omega,
\end{align}\end{subequations}
because of \ef{pomt},  $\bar\rho_0(y) =(\bar{A}- \bar{B}|{y}|^2 )^{{1}/({\ga-1})}$ and the Piola identity $\pl_k(JA^k_i)=0$. Here and thereafter, we use $\da_i^k $  to denote the Kronecker Delta symbol satisfying
$\da_i^k=1$ if $i=k$, and $\da_i^k=0$ if $i\neq k$.

We let $\pl_k=\frac{\pl}{\pl y_k}$,  $\pl^\al=\pl_1^{\al_1}\pl_2^{\al_2}\cdots\pl_n^{\al_n}$ for multi-index  $\al=(\al_1,\al_2, \cdots,\al_n)$ and $\pl^j=\sum_{|\al|=j}\pl^\al$ for nonnegative integer $j$.
We use $\bar \pl= y_1 \pl_{2}-y_2 \pl_{1}$ for  $n=2$,
and  $(\bar\pl_1,\bar\pl_2,\bar\pl_3)=y\times ( \pl_1, \pl_2, \pl_3)$ for $n=3$, respectively,
to denote  the angular momentum derivative, and let for $n=3$, similarly,
 $\bar\pl^\al=\bar\pl_1^{\al_1}\bar\pl_2^{\al_2}\bar\pl_3^{\al_3}$ for multi-index  $\al=(\al_1,\al_2, \al_3)$ and $\bar\pl^j=\sum_{|\al|=j}\bar\pl^\al$ for nonnegative integer $j$.
The divergence of a vector filed $F$ is ${\rm div} F=\da^{k}_i \pl_k F^i$,  the curl of a vector filed $F$ for $n=2$ is ${\rm curl} F=\pl_1 F_2- \pl_2 F_1$,
and the $i$-th component of the curl of a vector filed $F$ for $n=3$ is
$   [ {\rm curl} F ]_i = \epsilon^{ijk} \pl_j F_k$,
where  $\epsilon^{ijk}$ is the standard permutation symbol given by
\begin{align*}
\epsilon^{ijk}=
\begin{cases}
&1,    \ \  \ \  \textrm{even permutation of} \ \ \{1,2,3\},   \\
&-1,   \ \  \  \textrm{odd permutation of} \ \ \{1,2,3\}, \\
&0,    \ \ \  \  \textrm{otherwise}.
\end{cases}
\end{align*}
Along the flow map $\eta$, the $i$-th component of the gradient  of a function $f$ is
$ \lt[\nabla_\eta f \rt]_i = A^k_i \pl_k f$,
the divergence of a vector filed $F$ is ${\rm div}_\eta F = A^k_i \pl_k F^i$, the curl of a vector filed $F$ for $n=2$ is ${\rm curl}_\eta F=\lt[\nabla_\eta F_2 \rt]_1-\lt[\nabla_\eta F_1 \rt]_2=A^k_1\pl_k F_2-  A^k_2 \pl_k F_1$,
and the $i$-th component of the curl of a vector filed $F$ for $n=3$ is
 $
[ {\rm curl}_\eta F ]_{i}=\epsilon^{ijk} \lt[\nabla_\eta F_k \rt]_j   =\epsilon^{ijk} A^r_j \pl_r F_k$.

We introduce $\iota=(\ga-1)^{-1}$, $\sigma(y) =\bar\rho_0^{\ga-1}(y)=\bar A-\bar B |{y}|^2$, and
\begin{align*}
 &  \mathscr{E}^{m,i,j}(t)= (1+t)^{2m+(1-\da_{0m})\kappa} \lt\{ (1+t)^{1+\la} \lt\| \sa^{\frac{\iota+i}{2}}   \pl_t^{m+1} \pl^i \bar\pl^j \oa \rt\|_{L^2(\Omega)}^2 \rt.  \notag \\
 &  \ \lt. +   \lt\| \sa^{\frac{\iota+i}{2}} \pl_t ^m      \pl^i \bar\pl^j \oa \rt\|_{L^2(\Omega)}^2 +\lt\| \sa^{\frac{\iota+i+1}{2}}   \pl_t^m \pl^{i+1} \bar\pl^j \oa\rt\|_{L^2(\Omega)}^2  \rt\}
    \end{align*}
for nonnegative integers $m,i,j$; and
define the higher order weighted Sobolev norm by
$$
\mathscr{E}(t)=\sum_{0\le m+i+ j\le [\iota]+n+4} \mathscr{E}^{m,i,j}(t)+\sum_{ m+i+ j = [\iota]+n+4} \mathscr{E}^{m+1,i,j}(t).
$$
Here and thereafter, we use $\da_{ab}$ to denote the Kronecker Delta symbol satisfying $\da_{ab}=1$ if $a=b$, and $\da_{ab}=0$ if $a\neq b$.
We are  now ready to state the main result.

\begin{thm}\label{mainthm} Let $n=2,3$ and $0\le \la<1$.
There exists a positive constant $\bar \vea$ depending only on $n,\la,\ga$ and $M$ such that if $\mathscr{E}(0)\le \bar \vea$, then the problem \ef{system} admits a global  smooth solution in $[0,\iy)\times \Oa$ satisfying
\be\label{4.24-3}
\mathscr{E}(t)\le C \mathscr{E}(0) \ \ {\rm for} \ \  t \ge 0,
\ee
where $C$ is a positive constant which only depend on $n,\la,\ga$ and $M$ but do not depend on the time $t$.
\end{thm}

\begin{rmk}
When $m\ge 1$, the growth exponent of time in
$\mathscr{E}^{m,i,j}(t)$, $2m+\ka$, can be enhanced to $2m+2\ka -\vea$ for any fixed $\vea\in (0,\ka]$ by going through the entire proof, especially focusing on \ef{3.18-1}.
\end{rmk}

As a corollary of Theorem \ref{mainthm}, we have the following theorem for solutions to the original vacuum free boundary problem \ef{2.1}:

\begin{thm}\label{orig} Let $n=2,3$ and $0\le \la<1$.
There exists a positive constant $\bar \vea$ depending only on $n,\la,\ga$ and $M$ such that if $\mathscr{E}(0)\le \bar \vea$, then the problem \ef{2.1} admits a global  smooth solution $(\rho,u, \Oa(t))$ for
$t\in [0,\iy)$ satisfying that for $(t,y)\in [0,\iy)\times \Oa$,
 \begin{subequations}\label{4.23}
\begin{align}
& |x(t,y)- \bar x(t,y)|  \notag \\
& \quad  \le C (1+t)^{\ka}   \lt(  \sqrt{\mathscr{E}(0)}
 +  (1+t)^{\la -1}\lt(1+\da_{0\la} \ln(1+t) \rt)   \rt)
,\label{4.23-a}\\
&|\rho(t,x(t,y))-\bar\rho(t, \bar x (t,y))|  \notag \\
& \quad \le C
(1+t)^{-n\ka}  \bar\rho_0(y)  (  \sqrt{\mathscr{E}(0)}
+ (1+t)^{\la -1}\lt(1+\da_{0\la} \ln(1+t) \rt) ) ,\label{4.23-b}\\
&|u(t,x(t,y))-\bar u(t, \bar x (t,y))| \notag \\
& \quad \le C
(1+t)^{\ka-1}   \lt(  \sqrt{\mathscr{E}(0)}
 + (1+t)^{\la -1} \lt(1+\da_{0\la} \ln(1+t) \rt)  \rt)
, \label{4.23-c}
\end{align}
\end{subequations}
where $C$ is a positive constant independent of $t$ and $y$.
\end{thm}

\section{Proof of Theorem \ref{mainthm}}\label{sec3}
The proof of the global existence of smooth solutions is based on  the local existence theory (cf. \cite{10',16'}), together with
the following a priori estimates:
\begin{prop}\label{thm3.1}
Let $\oa(t,y)$ be a solution to  problem \ef{newsystem} in the time interval $[0, T]$ satisfying the following a priori assumptions:
\begin{align}
 \mathscr{E}(t) \le \ea_0^2 , \ \  t\in [0,T], \label{assume}
\end{align}
for some suitably small fixed positive number $\ea_0$ independent of $t$.
Then we obtain
\begin{align}\label{energy}
   \mathscr{E}  (t) \le C   \mathscr{E} (0), \ \  t\in [0,T],
\end{align}
where $C$ is a positive constant independent of $t$, which only depend on $n,\la,\ga$ and $M$.
\end{prop}
The proof of this proposition consists of curl and energy estimates, which are given in Lemmas
\ref{prop-curl} and \ref{3.31-1}, respectively.

\subsection{Preliminaries}
To simplify the presentation, we introduce some notation. Throughout the rest of paper,   $C$  will denote a positive constant which  only depend on the parameters of the problem, $n,\la,\ga$ and $M$, but does not depend on the data. They are referred as universal and can change
from one inequality to another one.
We will employ the notation $a\lesssim b$ to denote $a\le C b$, $a  \thicksim b$ to denote $C^{-1}b\le a\le Cb$, and $a\gtrsim b$ to denote $a\ge C^{-1} b$,
where  $C$ is the universal constant  as defined above.
Also we  will use $C(\beta)$ to denote  a certain positive constant
depending on quantity $\beta$.
We will use
$$ \int=\int_{\Omega},\ \ \|\cdot\|_{W^{k,q}}=\|\cdot\|_{W^{k,q}(\Omega)}, \ \  \|\cdot\|_{H^k}=\|\cdot\|_{W^{k,2}}   $$
for any constants $k\ge 0$ and $q\ge 1$.
Finally, in order to make the following statement more accurate, we bring in
\begin{align*}
&\mathscr{E}_I(t)=\sum_{0\le m+i+ j\le [\iota]+n+4} \mathscr{E}^{m,i,j}(t) \le \mathscr{E}(t), \\
& \mathscr{E}_{II}(t)=\sum_{0\le m+i+ j\le [\iota]+n+4} \mathscr{E}^{m+1,i,j}(t)  \le \mathscr{E}(t).
\end{align*}

We recall the following estimate obtained in (3.6) of  \cite{HZeng} and  (3.12) of \cite{HZeng2023}, whose proof is based on \ef{wsv} and \ef{hard}.
\begin{align}
 &  \sum_{m+
 2i+j\le 4} \lt\| \pl_t^m \pl^{i}\bar\pl^{j} \oa \rt\|_{L^\iy }^2  + \sum_{\substack{
 m+  2i+j=5 }} \lt(  \lt\| \pl_t^m  \pl^{i}\bar\pl^{j} \oa \rt\|_{H^{\frac{n+[\iota]-\iota}{2}}}^2
   + \lt\| \sa  \pl_t^m \pl^{i}\bar\pl^{j} \oa \rt\|_{L^\infty}^2\rt)   \notag \\
&  + \sum_{\substack{
6 \le m+ 2i+j \\
m+ i+j\le  [\iota]+n+3 }}  \lt\|\sa^{\frac{m+ 2i+j-4}{2}} \pl_t^m \pl^{i}\bar\pl^{j} \oa \rt\|_{L^\iy }^2
  \les  (1+t)^{-2m-(1-\da_{0m})\kappa} \mathscr{E}_I(t), \label{25.2.17}
\end{align}
provided that $\mathscr{E}_I(t)$ is finite.
This, together with \ef{assume}, means that
\bee
|\pl\oa(t,y)| \les \ea_0   \ \ {\rm for} \ \  (t,y)\in [0,T]\times \Oa.
\eee
Since $JA$ is the adjugate matrix of $(\frac{\pl \eta}{\pl y})$ and $\eta(t,y)=y+\oa(t,y)$, then
\bee
J A = \lt(\frac{\pl \eta}{\pl y} \rt)^*  =   \lt(1 +{\rm div}\oa\rt)\mathbb{I} -  \lt(\frac{\pl\oa}{\pl y}\rt) + \da_{n3} B,
\eee
where $B$  is the adjugate matrix of $(\frac{\pl \oa}{\pl y})$ given by
\bee
B= \lt(\frac{\pl\oa}{\pl y} \rt)^* = \lt[\begin{split} \pl_2 \omega  \times \pl_3 \omega \\
\pl_3 \omega  \times  \pl_1\omega  \\
\pl_1 \omega  \times \pl_2 \omega \end{split}\rt]  .
\eee
This, together with the fact that $(\frac{\pl x}{\pl y})(\frac{\pl x}{\pl y})^*=J \mathbb{I} $, implies that
 \begin{align}\label{Jaco}
 J=1+{\rm div} \oa + 2^{-1}\lt(|{\rm div} \oa|^2 + |{\rm curl} \oa|^2- |\pl \oa|^2\rt) + 3^{-1}  \da_{n3} B^s_r \pl_s \oa^r .
\end{align}
So, we have  that  for $t\in [0,T]$,
\be\label{7.9}
\|J-1\|_{L^\iy} \les \|\pl \oa\|_{L^\iy}  \les \ea_0  \ {\rm and}  \   \|A-  \mathbb{I}\|_{L^\iy} \les \|\pl \oa\|_{L^\iy}  \les \ea_0,
\ee
which gives, with the aid of the smallness of $\ea_0$, that  for $(t,y)\in [0,T]\times \Oa$,
 \be\label{6.7-1a}
  2^{-1}\le J \le 2     \ \  {\rm and} \ \
\max_{1\le i\le n, \
 1\le j\le n}  \lt| A^i_j \rt|\le 2 .
    \ee
 Moreover,  we have
$
   |[\na_\eta f]_i -\pl_i f| =|( A^r_i - \da^r_i) \pl_r f| \les \ea_0 |\pl f|
$ for any function $f$,
which means
 \be\label{6.7-1c}
   2^{-1} |\pl f| \le |\na_\eta f| \le 2 |\pl f|.
\ee

We recollect the identities in Lemma 4.4 of \cite{HZeng}, which indicate how the higher order functional are constructed.
For any vector field  $F $ with $F^i=F_i$, we have
\begin{subequations}\label{4.9-1}
\begin{align}
& A^k_rA^s_i (\pl_s  F^r)  \pl_t \pl_k F^i  =2^{-1} \pl_t \lt( |\nabla_\eta F|^2 - |{\rm curl}_\eta F|^2 \rt) +    \lt[ \nabla_{\eta } F^r \rt]_i \lt[\nabla_\eta \pl_t \oa^s\rt]_r\lt[\nabla_\eta F^i\rt]_s, \label{nabt}\\
&A^k_r A^s_i ( \pl_s F^r ) \pl_k F^i  = |\nabla_\eta F|^2 - |{\rm curl}_\eta F|^2. \label{nab}
\end{align}
\end{subequations}
We look back on the estimates for    commutators   in Lemmas 4.5-4.6 of \cite{HZeng}.  For any   function $f$ and multi-indexes  $\alpha$, $\alpha_1$ and $\beta$, we have
\begin{subequations}\label{4.9}\begin{align}
& \lt|\pl^{\alpha_1}\lt[\bar\pl^\beta , \pl^\alpha \rt]f\rt| \le C(\alpha_1,\alpha, \beta) \sum_{0\le j  \le |\beta|-1} \lt|\pl^{|\alpha_1|+|\alpha|} \bar\pl^j f \rt| , \label{commutator2}
\\
& \lt| \pl^\alpha \bar\pl^\beta \lt(  \sa^{-\iota}\pl_k  ( \sa^{\iota+1} f  ) \rt)- \sa^{-\iota - |\alpha|} \pl_k \lt(\sa^{\iota+ |\alpha|+ 1 }  \pl^\alpha \bar\pl^\beta f   \rt) \rt| \notag\\
&  \le   C \sum_{0\le j \le |\beta|-1} \lt(\sa \lt|\pl^{|\al|+1}\bar\pl^j f\rt| +   \lt|\pl^{|\al| }\bar\pl^j f \rt| \rt) +  C |\alpha| \sum_{0\le j \le |\beta|+1}
\lt| \pl^{|\alpha|-1} \bar\pl^j f \rt|,  \label{est1}
\end{align}\end{subequations}
where  $k=1,\cdots, n$ and $C=C(\al,\ba,\iota,\Omega)$ in \ef{est1}.
Indeed, the proof of \ef{4.9-1} and \ef{4.9} for $n=2$ is the same as that for $n=3$ shown in  \cite{HZeng}.

The differentiation formulae for $A$ and $J$ are
\begin{subequations}\label{7.12}\begin{align}
&\pl_j J=J A^s_r \pl_{j }\pl_{s } \oa^{r},     \   \bar\pl_j J=J A^s_r \bar\pl_j \pl_s  \oa^{r},  \      \pl_t J=J A^s_r   \pl_s v^r,  \label{7.12-1}\\
&\pl_j A^k_i = - A^k_r   A^s_i  \pl_{j }\pl_{s } \oa^r,  \       \bar\pl_j A^k_i = - A^k_r   A^s_i \bar\pl_{j }\pl_{s } \oa^r,   \
  \pl_t  A^k_i = - A^k_r A^s_i  \pl_s  v^r,\label{7.12-2}
\end{align}\end{subequations}
which, together with  the mathematical induction and \ef{6.7-1a}, implies that for  any polynomial function $\mathscr{P}$ and nonnegative integers $m,i,j$,
\begin{align}\label{25.2.21}
\lt|\pl_t^m  \pl^i \bar\pl^j  \mathscr{P} (A )  \rt| + \lt|\pl_t^m \pl^i \bar\pl^j  \mathscr{P}(J) \rt|
\les  \mathcal{I}^{m,i,j},
\end{align}
where $\mathcal{I}^{m,i,j}$ are defined inductively as follows:
\begin{subequations}\label{n5.30}\begin{align}
& \mathcal{I}^{0,0,0}=1, \\
& \mathcal{I}^{m,i,j} =
\sum_{\substack{0\le r\le m, \ 0\le k \le i, \ 0\le l \le j
\\ 0\le r+k+l \le m+i+j-1
  }} \mathcal{I}^{r,k,l} \lt|\pl_t^{m-r}\pl^{i-k}\bar\pl^{j-l}\pl\oa\rt|
. \label{5.30b}
 \end{align}\end{subequations}
Subsequently, we possess the following estimates which are necessary for deriving the curl and energy estimates.
\begin{lem} It holds that for any nonnegative integers $m,i,j$, and $t\in [0,T]$,
\begin{subequations}\label{25.3.1}\begin{align}
 & (1+t)^{2m+2+\kappa} \sum_{\substack{r\le m, \  k \le i, \  l \le j
\\ 1\le r+k+l
  }}\lt\|\sa^{\frac{\iota+i+1}{2}}
  \mathcal{I}^{r,k,l}
\pl_t^{m-r} \pl^{i-k}\bar\pl^{j-l}
\pl\pl_t\oa\rt\|_{L^2}^2 \notag\\
& \  \les    \mathscr{E}_I(t)
 \sum_{r+k+l\le m+i+j-1}\mathscr{E}^{r+1,k,l}(t),\ \ 1\le m+i+j \le [\iota]+n+4,\label{2.26-3}\\
   & (1+t)^{2m+2+\kappa} \sum_{\substack{r\le m, \  k \le i, \  l \le j
\\ 1\le r+k+l
  }}\lt\|\sa^{\frac{\iota+i}{2}}
  \mathcal{I}^{r,k,l}
\pl_t^{m-r} \pl^{i-k}\bar\pl^{j-l}
\pl\pl_t\oa\rt\|_{L^2}^2 \notag\\
& \  \les    \mathscr{E}_I(t)
 \sum_{r+k+l\le m+i+j}\mathscr{E}^{r+1,k,l}(t),\ \ 1\le m+i+j \le [\iota]+n+3,\label{2.26-4}\\
  & (1+t)^{2m+2+\kappa} \sum_{\substack{r\le m, \  k \le i, \  l \le j
\\ 1\le r+k+l
  }}\lt\|\sa^{\frac{\iota+i+1}{2}}
  \mathcal{I}^{r,k,l}
\pl_t^{m-r} \pl^{i-k}\bar\pl^{j-l}
\pl\pl_t \oa\rt\|_{L^2}^2 \notag\\
& \  \les    \mathscr{E}_I(t)
 \sum_{r+k+l\le m+i+j-1}\mathscr{E}^{r+1,k,l}(t),\ \   m+i+j = [\iota]+n+5,  \  1\le m,\label{4.9-3}
\end{align}\end{subequations}
provided that $\mathscr{E}(t)$ is small.
\end{lem}

{\it Proof}.  First,  \ef{2.26-3} can be  demonstrated  in spirit of  Lemma 3.4 of  \cite{HZeng2023}, that also includes the following estimate:
\begin{align}
  & (1+t)^{2m+(1-\da_{0m})\kappa} \sum_{\substack{r\le m, \  k \le i, \  l \le j
\\ 1\le r+k+l \le  m+ i+j-1
  }}\lt\|\sa^{\frac{\iota+i+1}{2}}
  \mathcal{I}^{r,k,l}
\pl_t^{m-r} \pl^{i-k}\bar\pl^{j-l}
\pl\oa\rt\|_{L^2}^2 \notag\\
& \ \  \les    \mathscr{E}_I(t)
 \sum_{r+k+l\le m+i+j-1}\mathscr{E}^{r,k,l}(t),\ \ \ \   2\le m+i+j \le [\iota]+n+5.\label{2.26-2}
\end{align}
Indeed, the time growth exponents of \ef{2.26-3} and  \ef{2.26-2} can be  easily observed, although they are not involved in Lemma 3.4 of  \cite{HZeng2023}.  The proof of  \ef{2.26-2} can  be found in Lemma 4.7 of \cite{HZeng} as well, where  time growth exponents are involved.

The proof of \ef{2.26-4} can be derived similarly as that of \ef{2.26-3}, with slight modifications involving \ef{wsv} and \ef{hard} in more places.

To prove \ef{4.9-3}, we notice that $\mathscr{E}_{I}(t)$ becomes
$\mathscr{E}_{II}(t)$ when $\oa$ is substituted by $\pl_t \oa$, regardless of the time growth exponents.  For example, we  replace  $\oa$ by $\pl_t \oa $  in \ef{25.2.17} to obtain
\begin{align}
 &     \sum_{m+
 2i+j\le 4} \lt\| \pl_t^m \pl^{i}\bar\pl^{j} \pl_t \oa \rt\|_{L^\iy }^2  + \sum_{\substack{
 m+  2i+j=5 }}  (  \lt\| \pl_t^m  \pl^{i}\bar\pl^{j}  \pl_t  \oa \rt\|_{H^{\frac{n+[\iota]-\iota}{2}}}^2
 \notag \\
&  \ \  + \lt\| \sa  \pl_t^m \pl^{i}\bar\pl^{j}  \pl_t  \oa \rt\|_{L^\infty}^2 )     + \sum_{\substack{
6 \le m+ 2i+j \\
m+ i+j\le  [\iota]+n+3 }}  \lt\|\sa^{\frac{m+ 2i+j-4}{2}} \pl_t^m \pl^{i}\bar\pl^{j}  \pl_t  \oa \rt\|_{L^\iy }^2
\notag \\
& \ \
  \les  (1+t)^{-2m-2-\kappa} \mathscr{E}_{II}(t).
  \label{4.9-2}
\end{align}
Then,  we utilize the identical method of demonstration as that employed for \ef{2.26-2} to establish
\begin{align}
  & (1+t)^{2m+2+\kappa} \sum_{\substack{r\le m, \  k \le i, \  l \le j
\\ 1\le r+k+l \le  m+ i+j-1
  }}\lt\|\sa^{\frac{\iota+i+1}{2}}
  \mathcal{I}^{r,k,l}
\pl_t^{m-r} \pl^{i-k}\bar\pl^{j-l}
\pl\pl_t\oa\rt\|_{L^2}^2 \notag\\
& \ \  \les    \mathscr{E}_I(t)
 \sum_{r+k+l\le m+i+j-1}\mathscr{E}^{r+1,k,l}(t),\ \ \ \   2\le m+i+j \le [\iota]+n+5.\label{4.9-4}
\end{align}
When $m+i+j = [\iota]+n+5$ and $m \ge 1$, it follows from \ef{commutator2} and \ef{2.26-2} that
\begin{align}
  & \lt\|\sa^{\frac{\iota+i+1}{2}}
  \mathcal{I}^{m,i,j}
\rt\|_{L^2}^2
\les    \lt\|\sa^{\frac{\iota+i+1}{2}}\pl_t^{m}\pl^{i}\bar\pl^{j}\pl\oa
\rt\|_{L^2}^2\notag \\
& + \sum_{\substack{  r\le m, \  k \le i, \   l \le j
\\ 1\le r+k+l \le m+i+j-1
  }} \lt\|\sa^{\frac{\iota+i+1}{2}} \mathcal{I}^{r,k,l} \pl_t^{m-r}\pl^{i-k}\bar\pl^{j-l}\pl\oa\rt\|_{L^2}^2\notag \\
&  \les (1+t)^{-2m-\ka }\{\sum_{l\le j} \mathscr{E}^{m,i,l} (t)
+ \mathscr{E}_I (t)
 \sum_{r+k+l\le [\iota]+n+4}\mathscr{E}^{r,k,l}(t)  \}, \notag
\end{align}
which, together with \ef{25.2.17} and \ef{4.9-2}, implies that
\begin{align}
  & \lt\|\sa^{\frac{\iota+i+1}{2}}
  \mathcal{I}^{m,i,j}
\pl\pl_t\oa\rt\|_{L^2}^2
\les \lt\|\sa^{\frac{\iota+i+1}{2}}
  \mathcal{I}^{m,i,j} \rt\|_{L^2}^2 \lt\|
\pl\pl_t\oa\rt\|_{L^\iy}^2
 \notag\\
&    \les   (1+t)^{-2m-2-2\ka }\{\mathscr{E}_I(t) \sum_{l\le j} \mathscr{E}^{m,i,l}(t)+\mathscr{E}_I^2 (t) \mathscr{E}_{II} (t)\}
\notag\\
& \les  (1+t)^{-2m-2-2\ka }\mathscr{E}_I(t)(1+\mathscr{E}_I(t))
 \sum_{r+k+l\le m+i+j-1}\mathscr{E}^{r+1,k,l}(t). \notag
\end{align}
When this is taken with \ef{4.9-4} and the smallness of $\mathscr{E}_I(t)$, it proves \ef{4.9-3}.
\hfill$\Box$

\subsection{Curl  estimates}\label{sec-curl}
\begin{lem}\label{prop-curl}
Let $m,i,j$ be nonnegative integers, then it holds that for $m+ i+j\le [\iota]+n+4$ and $t\in [0,T]$,
\begin{align}
&(1+t)^{2m+4+\ka-2\la} \lt\|\sa^{\frac{\iota+i+1}{2}}  \pl_t^{m}\pl^{i}\bar\pl^{j}{\rm curl}_\eta \pl_t \oa \rt\|_{L^2}^2
(t)
\notag\\
&  \les\lt\|\sa^{\frac{\iota+i+1}{2}}  \pl^{i}\bar\pl^{j}{\rm curl}_\eta  \pl_t \oa \rt\|_{L^2}^2
(0)
+\ea_0^2 \sum_{r+k+l \le m+i+j }  \sup_{\tau\in [0,t]} \mathscr{E}^{r+1,k,l} (\tau) . \label{curl-a}
\end{align}
\end{lem}

{\em Proof}.  We use  the Piola identity $\pl_k(JA^k_i)=0$ and $\bar\rho_0(y) =(\bar{A}- \bar{B}|{y}|^2 )^{{1}/({\ga-1})}$ to express  equation \ef{3-1-3}   as
 \bee\begin{split}
&\ta   \pl_t^2 \oa + \lt ( {(1+t)^{-\la}}\ta+ 2\ta_t \rt)  \pl_t \oa    + \ka   \ta^{n-1-n\ga}  \eta
 \\
 &+  (1+\iota)  \ta^{n-1-n\ga}  \nabla_\eta   \lt(\bar\rho_0^{\ga-1} {J}^{1-\ga}\rt)  = 0.
\end{split}\eee
Let ${\rm curl}_\eta$ act on it, and use the fact ${\rm curl}_\eta \eta =0$ and ${\rm curl}_\eta \nabla_\eta =0$ to give
$$
\ta {\rm curl}_\eta \pl_t^2 \oa  +\lt ( {(1+t)^{-\la}}\ta+ 2\ta_t \rt)  {\rm curl}_\eta   \pl_t \oa    =0.
$$
Commuting $\pl_t$ with ${\rm curl}_\eta$ and noting the integrating-factor, we have
 \begin{align}
&{\rm curl}_\eta   \pl_t \oa  =\theta^{-2}(t) e^{-\frac{1}{1-\la}(1+t)^{1-\la}}\lt\{
e^{\frac{1}{1-\la}}{\rm curl}_{\eta}     \pl_t \oa \big|_{t=0}\rt.\notag\\
&\   \lt. +\int_0^t    \ta^2(\tau) e^{\frac{1}{1-\la}(1+\tau)^{1-\la}} \lt[\pl_\tau , {\rm curl}_\eta  \rt]    \pl_\tau \oa  d\tau
\rt\}.\label{3-5-4-a}
 \end{align}

To prove \ef{curl-a} for $m=0$, we
take $\pl^{\al}\bar\pl^{\ba}$  onto  \ef{3-5-4-a} to obtain
\begin{align}
& \pl^{\al}\bar\pl^{\ba} {\rm curl}_\eta    \pl_t \oa =  \pl^{\al}\bar\pl^{\ba} {\rm curl}_{\eta}     \pl_t \oa \big|_{t=0} \theta^{-2}(t) e^{-\frac{1}{1-\la}\lt((1+t)^{1-\la}-1\rt)}
\notag\\
&   +\theta^{-2}(t) e^{-\frac{1}{1-\la}(1+t)^{1-\la}} \int_0^t    \ta^2(\tau) e^{\frac{1}{1-\la}(1+\tau)^{1-\la}} \pl^{\al}\bar\pl^{\ba} \lt[\pl_\tau , {\rm curl}_\eta  \rt]   \pl_\tau \oa  d\tau,
  \label{8.4c}
\end{align}
where $\pl^{\al}\bar\pl^{\ba} \lt[\pl_\tau , {\rm curl}_\eta  \rt]   \pl_\tau \oa$ can be  written  as
$$
 \pl^{ \al }\bar\pl^{ \ba } \lt(  \pl_\tau [{\rm curl}_\eta \pl_\tau \oa]_l  -   [{\rm curl}_\eta \pl_\tau^2 \oa ]_l   \rt)
=\pl^{ \al }\bar\pl^{ \ba }(\epsilon^{ljk} (\pl_\tau  A^r_j)
\pl_\tau  \pl_r \oa_k).
$$
This, together with \ef{4.9-3}, implies  that for $h+|\al|+|\ba|\le [\iota]+n+4$,
\begin{align}
&  \lt\| \sa^{\frac{\iota+|\al|+1}{2}}\pl_\tau^h \pl^{\al}\bar\pl^{\ba} \lt[\pl_\tau , {\rm curl}_\eta  \rt]   \pl_\tau \oa  \rt\|_{L^2}^2
\notag\\
&  \les   \sum_{r\le h, \ k\le |\al|,\ l \le |\ba|} \lt\|\sa^{\frac{\iota+|\al|+1}{2}}
 \mathcal{I}^{1+r,k,l}
\pl_\tau^{h-r} \pl^{|\al|-k}\bar\pl^{|\ba|-l}
\pl \pl_\tau \oa \rt\|_{L^2}^2
\notag\\
& \le \sum_{1\le r\le h+1
, \ k\le |\al|,\ l \le |\ba|} \lt\|\sa^{\frac{\iota+|\al|+1}{2}}
 \mathcal{I}^{r,k,l}
\pl_\tau^{h+1-r} \pl^{|\al|-k}\bar\pl^{|\ba|-l}
\pl \pl_\tau \oa \rt\|_{L^2}^2
\notag\\
& \les (1+\tau)^{-2h-4-\ka}  \ea_0^2  \sum_{r+k+l \le h+|\al|+|\ba|  } \mathscr{E}^{r+1,k,l} (\tau) . \label{3.2-1}
\end{align}
It should be pointed  that for any fixed $k >0$,
\begin{align}
&  e^{-\frac{1}{1-\la}(1+t)^{1-\la}} \int_0^t     e^{\frac{1}{1-\la}(1+\tau)^{1-\la}}
(1+\tau)^{-k}     d\tau
\notag\\
& \le
e^{-\frac{1}{1-\la}\lt(
(1+t)^{1-\la}
-\lt(1+\frac{t}{2}\rt)^{1-\la} \rt)}
\int_0^{t/2} (1+\tau)^{-k} d\tau
\notag\\
& \ +e^{-\frac{1}{1-\la}(1+t)^{1-\la}}
(1+t)^{\la}\lt(1+\frac{t}{2}\rt)^{-k} \int_{t/2}^{t} d e^{\frac{1}{1-\la}(1+\tau)^{1-\la}}
\notag\\
& \le C(\la, k) \lt(1+t\rt)^{\la -k}, \label{25.2.12}
\end{align}
which, along with  \ef{8.4c},  \ef{3.2-1} and \ef{5.26}, proves
   \ef{curl-a} for $m=0$.

When $m\ge1$, we  apply $\pl_t^{m}$ onto \ef{8.4c} to obtain
\begin{align}
& \pl_t^{m}\pl^{\al}\bar\pl^{\ba} {\rm curl}_\eta    \pl_t \oa  =  \pl^{\al}\bar\pl^{\ba} {\rm curl}_{\eta}   \pl_t \oa   \big|_{t=0} \pl_t^m \lt(\theta^{-2}(t) e^{-\frac{1}{1-\la}\lt((1+t)^{1-\la}-1\rt)} \rt)
\notag\\
&+ \sum_{0\le h \le m}  \frac{m!}{h!(m-h)!}    \pl_t^{m-h }\theta^{-2}(t)  \pl_t^h \lt( e^{-\frac{1}{1-\la}(1+t)^{1-\la}}
\rt. \notag\\
& \lt. \times  \int_0^t   e^{\frac{1}{1-\la}(1+\tau)^{1-\la}}  \ta^2(\tau) \pl^{\al}\bar\pl^{\ba} \lt[\pl_\tau , {\rm curl}_\eta  \rt]   \pl_\tau \oa    d\tau\rt) .
  \label{2.28}
\end{align}
It follows from the integration by parts that for any function $f(t,y)$,
\begin{align}
&\pl_t   \lt( e^{-\frac{1}{1-\la}(1+t)^{1-\la}} \int_0^t    e^{\frac{1}{1-\la}(1+\tau)^{1-\la}} f(\tau,y)   d\tau \rt)
\notag \\
&=f(t,y)- (1+t)^{-\la}  e^{-\frac{1}{1-\la}(1+t)^{1-\la}} \int_0^t    e^{\frac{1}{1-\la}(1+\tau)^{1-\la}} f(\tau,y)   d\tau \notag\\
&=f(t,y)- (1+t)^{-\la}  e^{-\frac{1}{1-\la}(1+t)^{1-\la}} \int_0^t  (1+\tau\emph{})^{\la}    f(\tau,y)   d e^{\frac{1}{1-\la}(1+\tau)^{1-\la}}\notag\\
&= (1+t)^{-\la}  e^{-\frac{1}{1-\la} \lt((1+t)^{1-\la}-1\rt)} f(0,y)    \notag \\
&    +(1+t)^{-\la}  e^{-\frac{1}{1-\la}(1+t)^{1-\la}} \int_0^t  e^{\frac{1}{1-\la}(1+\tau)^{1-\la}} \pl_\tau\lt( (1+\tau)^{\la}    f(\tau,y) \rt)  d\tau, \notag
\end{align}
which proves  \ef{curl-a} for $m\ge1$, using \ef{3.2-1}-\ef{2.28} and \ef{5.26}.
\hfill $\Box$

\subsection{Energy estimates}
\begin{lem}\label{3.31-1}
Let $ m\ge 1,i,j$ be nonnegative integers, then
it holds that for  $m+i+j\le [\iota]+n+5$  and $t\in [0,T]$,
\begin{align}\label{3.28-2}
&\mathscr{E}^{m,i,j}(t)
+\int_0^t  \{ (1+\tau)^{2m+\ka +1} \int  \sigma^{  \iota + i}
 |  \pl_\tau^{m+1}\pl^i \bar\pl^j \oa  |^2 dy
\notag\\
& +(1+\tau)^{-1} \mathscr{E}^{m,i,j}(\tau)  \}d\tau
\les \sum_{r+k+l\le m+i+j-1} \mathscr{E}^{r+1,k,l}(0);
\end{align}
and for $i+j\le [\iota]+n+4$  and $t\in [0,T]$,
\be\label{3.31}
 \mathscr{E}^{0,i,j}(t)\les \mathscr{E}^{0,i,j}(0)+\sum_{r+k+l\le i+j} \mathscr{E}^{r+1,k,l}(0) .
\ee
\end{lem}

{\em Proof}.  Once \ef{3.28-2} is proved, \ef{3.31} follows easily from
$
\oa(t,y)=\oa(0,y)+\int_0^t \pl_\tau \oa (\tau, y) d\tau
$
and
\begin{align*}
& \lt\| \sa^{\frac{\iota+i}{2}}     \pl^i \bar\pl^j \oa \rt\|_{L^2} +\lt\| \sa^{\frac{\iota+i+1}{2}}  \pl^{i+1} \bar\pl^j \oa\rt\|_{L^2}
\\
& \les    \sqrt{ \mathscr{E}^{0,i,j}(0) }
 +   \int_0^t (1+\tau )^{-1-2^{-1}\ka} \sqrt{ \mathscr{E}^{1,i,j}(\tau) } d\tau \\
& \les  \sqrt{ \mathscr{E}^{0,i,j}(0) }+
 \sqrt{  \sum_{r+k+l\le i+j} \mathscr{E}^{r+1,k,l}(0) }.
    \end{align*}
Hence, we will focus on proving \ef{3.28-2} in what follows.
For this purpose, we apply $\pl_t^m\pl^\alpha \bar\pl^\beta$
onto the product of  $\sa^{-\iota}$ and \ef{3-1-3} to  obtain
\begin{align}
&\ta^{n\ga-n+2}     \pl_t^{m+2} \pl^\alpha \bar\pl^\beta \oa_i +  (  1+t  )  \pl_t^{m+1}  \pl^\alpha \bar\pl^\beta \oa_i +  (\ka+m  )\pl_t^m  \pl^\alpha \bar\pl^\beta \oa_i
 \notag \\
  &
   + \sigma^{ - \iota - |\alpha|} \pl_k ( \sigma^{\iota+ |\alpha|+1} ( \mathcal{R}_{1, i}^{m, \al,\ba,k}
    -{J}^{1-\ga}  (   A^k_r  A^s_i  \pl_s  \pl_t^m\pl^\alpha \bar\pl^\beta   \oa^r
 \notag\\&
 + \iota^{-1} A_i^k {\rm div}_\eta   \pl_t^m\pl^\alpha \bar\pl^\beta  \oa)))
  =-  (\ta^{n\ga-n+1}  ( (1+t)^{-\la}\ta+ 2\ta_t )  \notag\\
  &  - (1+t) +m  (\ta^{n\ga-n+2})_t  ) \pl_t^{m+1} \pl^\alpha \bar\pl^\beta \oa_i
+  \mathcal{R}_{2,i}^{m, \al,\ba}+\mathcal{R}_{3,i}^{m, \al,\ba}, \label{25.2.19}
\end{align}
where
\begin{align*}
&\mathcal{R}_{1,i}^{m,\al,\ba,k} =       \pl_t^m\pl^\alpha\bar\pl^\beta ( A_i^k {J}^{1-\ga}  -  \da_i^k ) \\
& \ + {J}^{1-\ga} (   A^k_r  A^s_i  \pl_s \pl_t^{m} \pl^\alpha \bar\pl^\beta   \omega^r   + \iota^{-1} A_i^k {\rm div}_\eta   \pl_t^{m} \pl^\alpha \bar\pl^\beta   \omega    ),
 \\
& \mathcal{R}_{2,i}^{m,\al,\ba}=   { \sigma^{ - \iota - |\alpha|}} \pl_k \lt( \sigma^{\iota+ |\alpha|+1}   \pl_t^m\pl^\alpha\bar\pl^\beta ( A_i^k {J}^{1-\ga}  -  \da_i^k  ) \rt)
\\
  &   \  -\pl_t^m\pl^\alpha \bar\pl^\beta \lt( \sigma^{-\iota } \pl_k  \lt( \sigma^{\iota+1}  (A_i^k {J}^{1-\ga} -  \da_i^k ) \rt) \rt),\\
   & \mathcal{R}_{3,i}^{m,\al,\ba}=-  \sum_{1\le h \le m} \frac{m!}{h!(m-h)!} ( \pl_t^{m-h+1} \pl^\alpha \bar\pl^\beta \oa_i )    \frac{d^h}{dt^h}  (\ta^{n\ga-n+1}  ( (1+t)^{-\la}\ta \\
   &    \ + 2\ta_t )    - (1+t)  )  -
   \sum_{2\le h \le m} \frac{m!}{h!(m-h)!}    (\pl_t^{m-h+2} \pl^\alpha \bar\pl^\beta \oa_i )
   \frac{d^h}{dt^h}\ta^{n\ga-n+2}.
\end{align*}
Based on the equation above, the proof of \ef{3.28-2} consists of the following three steps.

{\it Step 1}. Let  $b$ be a time decay rate, which is a nonnegative constant and will be determined later, then we integrate the product of $(1+t)^b \sigma^{  \iota + |\alpha|}   \pl_t^{m+1}  \pl^\alpha \bar\pl^\beta \oa^i$ and \ef{25.2.19}   over $\Omega$  and employ \ef{nabt}  to  achieve
\begin{align}\label{5.2}
\frac{d}{dt} \sum_{1\le j\le 3}\mathcal{E}_j^{m,\al,\ba}(t)+  \mathcal{D}_1^{m,\al,\ba}(t)  = \sum_{j=1,2} \frac{d}{dt} \mathcal{F}_j^{m,\al,\ba}(t) + \mathcal{H}_1^{m,\al,\ba}(t) ,
\end{align}
where
\begin{align}
&  \mathcal{E}_1^{m, \al,\ba}(t)   =    \frac{1}{2}(1+t)^b  \ta^{n\ga-n+2} \int \sigma^{  \iota + |\alpha|}
\lt|  \pl_t^{m+1}\pl^\alpha \bar\pl^\beta \oa \rt|^2 dy ,
\notag\\
&  \mathcal{E}_2^{m, \al,\ba}(t)   =  \frac{1 }{2}  (\ka+m ) (1+t)^b
 \int  \sigma^{  \iota + |\alpha|}\lt| \pl_t^{m} \pl^\alpha \bar\pl^\beta \oa\rt|^2 dy,
\notag \\
& \mathcal{E}_3^{m, \al,\ba}(t)   =\frac{1}{2}
 (1+t)^b \int  \sigma^{\iota+ |\alpha|+1} {J}^{1-\ga}  (   |\na_\eta \pl_t^m  \pl^\alpha \bar\pl^\beta \oa|^2   + \iota^{-1}    |{\rm div}_\eta  \pl_t^m  \pl^\alpha \bar\pl^\beta \oa|^2      ) dy
,
\notag\\
&  \mathcal{D}_1^{m,\al,\ba} (t)  =     (  1+t)^{b+1} \int  \sigma^{  \iota + |\alpha|}
\lt|  \pl_t^{m+1}\pl^\alpha \bar\pl^\beta \oa \rt|^2 dy  ,
\notag\\
&  \mathcal{F}_1^{m, \al,\ba}(t)   =    \frac{1}{2} (1+t)^b
  \int  \sigma^{\iota+ |\alpha|+1} {J}^{1-\ga}  |{\rm curl}_\eta \pl_t^m  \pl^\alpha \bar\pl^\beta \oa|^2
  dy,
\notag\\
&  \mathcal{F}_2^{m, \al,\ba}(t)   =  (1+t)^b   \int   \sigma^{\iota+ |\alpha|+1} \mathcal{R}_{1,i}^{m, \al,\ba,k}  \pl_k \pl_t^m  \pl^\alpha \bar\pl^\beta \oa^i dy,
\notag\\
& \mathcal{H}_1^{m,\al,\ba}(t)   =   - \int  \sigma^{\iota+ |\alpha|+1}   \pl_t (  (1+t)^b\mathcal{R}_{1,i}^{m,\al,\ba,k}) \pl_k \pl_t^m \pl^\alpha \bar\pl^\beta \oa^i dy
\notag\\
& \ +(1+t)^b \sum_{j=2,3} \int  \sigma^{  \iota + |\alpha|}  \mathcal{R}_{j,i}^{m, \al,\ba}  \pl_t^{m+1}  \pl^\alpha \bar\pl^\beta \oa^i dy
\notag\\
&\ +
\{\frac{1}{2} ((1+t)^b\ta^{n\ga-n+2})_t -(1+t)^b (\ta^{n\ga-n+1}  ( (1+t)^{-\la}\ta+ 2\ta_t )   \notag\\
&  \ - (1+t)
  +m  (\ta^{n\ga-n+2})_t  )  \}
\int \sigma^{  \iota + |\alpha|}
\lt|  \pl_t^{m+1}\pl^\alpha \bar\pl^\beta \oa \rt|^2 dy
\notag\\
&\
 + \frac{1 }{2}  (\ka+m ) ((1+t)^b)_t
 \int  \sigma^{  \iota + |\alpha|}\lt| \pl_t^{m} \pl^\alpha \bar\pl^\beta \oa\rt|^2 dy
 \notag\\
&\
+ \frac{1}{2}    \int  \sigma^{\iota+ |\alpha|+1}  \pl_t ((1+t)^b{J}^{1-\ga} )  (   |\na_\eta  \pl_t^m  \pl^\alpha \bar\pl^\beta \oa|^2 - |{\rm curl}_\eta \pl_t^m  \pl^\alpha \bar\pl^\beta \oa|^2
\notag\\
& \ + \iota^{-1}    |{\rm div}_\eta  \pl_t^m  \pl^\alpha \bar\pl^\beta \oa|^2     ) dy
  - (1+t)^b\int  \sigma^{\iota+ |\alpha|+1}  {J}^{1-\ga}  (
[\na_\eta  \pl_t^m  \pl^\alpha \bar\pl^\beta \oa^r]_i[\na_\eta   \pl_t\oa^s]_r
\notag\\
&\ \times [\na_\eta   \pl_t^m  \pl^\alpha \bar\pl^\beta \oa^i]_s
 -\iota^{-1} (\pl_t A_i^k) ( {\rm div}_\eta   \pl_t^m  \pl^\alpha \bar\pl^\beta \oa) \pl_k \pl_t^m  \pl^\alpha \bar\pl^\beta \oa^i ) dy.\notag
\end{align}
Clearly, the terms on the left-hand side of equation \ef{5.2} are the desired ones. As we proceed, we will address the terms on the right-hand side. Throughout this process, we will make frequent use of the estimates \ef{6.7-1a},  \ef{6.7-1c} and \ef{4.9}-\ef{25.3.1}.

Regarding $\mathcal{F}_1^{m, \al,\ba}$,  we  deduce from
 \ef{6.7-1a} and \ef{25.2.21}  that for  $m \ge 1$,
\begin{align*}
& \lt| {\rm curl}_\eta   \pl_t^{m } \pl^{\al}\bar\pl^{\ba}\oa -\pl_t^{m-1} \pl^{\al}\bar\pl^{\ba}{\rm curl}_\eta \pl_t\oa \rt|   \\
& \les  |\pl_t^m \pl^{ \al }[\pl, \bar\pl^{ \ba }]\oa|
+\sum_{\substack{r\le m-1, \ k\le |\al|, \ l \le |\ba|
\\ 1\le r+k+l
  }}  \mathcal{I}^{r,k,l}
\lt|\pl_t^{m-r}\pl^{|\al|-k}\bar\pl^{|\ba|-l}\pl \oa\rt| ,
\end{align*}
which, together with \ef{commutator2} and \ef{2.26-3}, implies that for  $m \ge 1$,
 \begin{align}
& (1+t)^{2m+\kappa} \int \sa^{ \iota+|\al|+1 } \lt| {\rm curl}_\eta   \pl_t^{m } \pl^{\al}\bar\pl^{\ba}\oa -\pl_t^{m-1} \pl^{\al}\bar\pl^{\ba}{\rm curl}_\eta \pl_t\oa   \rt|^2 dy \notag \\
&
\les \sum_{0\le j\le |\ba|-1} \mathscr{E}^{m,|\al|,j}(t) +  \ea_0^2 \sum_{r+k+l\le m+|\al|+|\ba|-2}\mathscr{E}^{r+1,k,l}(t). \label{8.8b}
\end{align}
Then, it follows from \ef{6.7-1a} and  \ef{curl-a}   that   for $m \ge 1$,
\begin{align}
&\mathcal{F}_1^{m, \al,\ba}(t)  \les    (1+t)^b
  \int  \sigma^{\iota+ |\alpha|+1}    |  \pl_t^{m-1} \pl^{\al}\bar\pl^{\ba}{\rm curl}_\eta \pl_t\oa|^2
  dy\notag\\
  & \ +
  (1+t)^b \int  \sigma^{\iota+ |\alpha|+1}    |{\rm curl}_\eta \pl_t^m  \pl^\alpha \bar\pl^\beta \oa- \pl_t^{m-1} \pl^{\al}\bar\pl^{\ba}{\rm curl}_\eta \pl_t\oa|^2
  dy \notag\\
 &  \les   (1+t)^{b-2m-\ka} \{ (1+t)^{2\la-2} (   \sum_{r+k+l \le m+|\al|+|\ba|-1 }  \sup_{\tau\in [0,t]} \ea_0^2   \mathscr{E}^{r+1,k,l} (\tau)  \notag\\
 &\ + \sum_{r+k+l\le |\al|+|\ba| }\mathscr{E}^{r+1,k,l}(0)  )  +
  \sum_{r+k+l\le m+|\al|+|\ba|-2}\mathscr{E}^{r+1,k,l}(t) \},
  \label{3.15-1}
 \end{align}
 where  \ef{commutator2} and \ef{2.26-3} have been used to rewrite the initial data in \ef{curl-a}.

For $\mathcal{F}_2^{m, \al,\ba}$, we  utilize \ef{6.7-1a}, \ef{commutator2}, \ef{7.12}, \ef{25.2.21} and
\begin{align*}
& \mathcal{R}_{1,i}^{m,\al,\ba,k} =       \pl_t^m\pl^\alpha\bar\pl^\beta ( A_i^k {J}^{1-\ga}  -  \da_i^k )  + {J}^{1-\ga} (   A^k_r  A^s_i   \pl_t^{m} \pl^\alpha \bar\pl^\beta  \pl_s \omega^r  \\
& \  + \iota^{-1} A_i^k A^s_r    \pl_t^{m} \pl^\alpha \bar\pl^\beta   \pl_s \omega^r    )     + {J}^{1-\ga} (   A^k_r  A^s_i \pl_t^{m} \pl^\alpha  [\pl_s,  \bar\pl^\beta ]  \omega^r  \\
& \ + \iota^{-1} A_i^k A^s_r    \pl_t^{m} \pl^\alpha [ \pl_s, \bar\pl^\beta ]   \omega^r    )
\end{align*}
to obtain that for $m\ge 1$,
\begin{align}
&|\mathcal{R}_{1,i}^{m,\al,\ba,k} | \les   \sum_{\substack{r\le m-1,\   k \le |\al|, \   l \le |\ba|
\\ 1\le r+ k+l
  }} \mathcal{I}^{r,k,l} \lt| \pl_t^{m-r} \pl^{|\al|-k}\bar\pl^{|\ba|-l}\pl \oa\rt|
  \notag \\
  &\
  + |\beta |  \sum_{  j  \le |\ba|-1}|\pl_t^{m} \pl^{|\alpha|+1} \bar\pl^{j} \oa|. \label{3.10-1}
 \end{align}
This, along with \ef{2.26-3},   gives rise to the inference that for $m\ge 1$,
\begin{align}
&(1+t)^{2m+\ka} \int \sigma^{  \iota + |\alpha|+1} |\mathcal{R}_{1,i}^{m,\al,\ba,k} |^2 dy \notag \\
   & \ \les
 \sum_{r+k+l\le m+|\al|+|\ba|-2} \ea_0^2 \mathscr{E}^{r+1,k,l}(t)  + |\beta |   \sum_{  j  \le |\ba|-1} \mathscr{E}^{m,|\al|,j}(t).\label{3.10-2}
 \end{align}
As a result, we have  that for any $\vea>0$,
\begin{align}\label{3.15-2}
 | \mathcal{F}_2^{m, \al,\ba} (t) |   \les \varepsilon \mathcal{E}_3^{m, \al,\ba}(t)  +\varepsilon^{-1}
(1+t)^{b-2m-\ka}
 \sum_{r+k+l\le m+|\al|+|\ba|-2}\mathscr{E}^{r+1,k,l} (t).
\end{align}

To deal with $\mathcal{H}_1^{m,\al,\ba}(t)$, we need to exercise more careful judgment in our estimates.
Just as in the process of derving \ef{3.10-1} and \ef{3.10-2}, we get, using \ef{4.9-3} additionally, that
\begin{align*}
&
|\pl_t \mathcal{R}_{1,i}^{\al,\ba,k} | \les     \sum_{\substack{r\le m ,\   k \le |\al|, \   l \le |\ba|
\\ 1\le r+ k+l
  }} \mathcal{I}^{r,k,l} \lt| \pl_t^{m+1-r} \pl^{|\al|-k}\bar\pl^{|\ba|-l}\pl \oa\rt| \\
   & \ +|\beta | \sum_{ j \le |\ba|-1}|\pl_t^{m+1} \pl^{|\alpha|+1} \bar\pl^{j} \oa|+|\pl v||\pl_t^{m} \pl^{|\alpha|+1} \bar\pl^{|\ba|} \oa|
\end{align*}
and
\begin{align*}
&(1+t)^{2m+2+\ka} \int \sigma^{  \iota + |\alpha|+1} |\pl_t \mathcal{R}_{1,i}^{m,\al,\ba,k} |^2 dy \\
   & \ \les
 \sum_{r+k+l\le m+|\al|+|\ba|-1}  \ea_0^2 \mathscr{E}^{r+1,k,l}(t)  + |\beta |   \sum_{  j  \le |\ba|-1}  \mathscr{E}^{m+1,|\al|,j}(t),
 \end{align*}
which, in conjunction with \ef{3.10-2}, suggests
 that for any $\vea>0$,
 \begin{align}
&\int  \sigma^{\iota+ |\alpha|+1} \big|  \pl_t (  (1+t)^b\mathcal{R}_{1,i}^{m,\al,\ba,k}) \pl_k \pl_t^m \pl^\alpha \bar\pl^\beta \oa^i \big|dy  \notag\\
& \les  \vea  (1+t)^{-1} \mathscr{E}^{m ,|\al|,|\ba|}(t)
 +\vea^{-1}(1+t)^{2b-4m-2\ka-1} \{|\beta |   \mathscr{E}^{m+1,|\al|,|\ba|-1}
\notag\\
 & \  +
 \sum_{r+k+l = m+|\al|+|\ba|-1}\ea_0^2\mathscr{E}^{r+1,k,l} +
 \sum_{r+k+l\le m+|\al|+|\ba|-2}\mathscr{E}^{r+1,k,l} \}(t).
 \label{3.14-2}
 \end{align}

Concerning the term that incorporates $\mathcal{R}_{2,i}^{m,\al,\ba}$, we use \ef{6.7-1a}, \ef{est1},  \ef{7.12} and \ef{25.2.21} to produce
\begin{align}
&  |\mathcal{R}_{2,i}^{m,\al,\ba} |
   \les |\ba|  \sum_{  j\le |\ba|-1}   \sum_{\substack{r\le m-1,\   k \le |\al|+1, \   l \le j
  }}  \sa  \mathcal{I}^{r,k,l} \lt| \pl_t^{m-r} \pl^{|\al|+1-k}\bar\pl^{j-l}\pl \oa\rt|
  \notag \\
&\   + |\ba|  \sum_{  j\le |\ba|-1} \sum_{\substack{r\le m-1,\   k \le |\al|, \   l \le j
  }}   \mathcal{I}^{r,k,l} \lt| \pl_t^{m-r} \pl^{|\al| -k}\bar\pl^{j-l}\pl \oa\rt| \notag \\
&\    +  |\al| \sum_{  j\le |\ba|+1} \sum_{\substack{r\le m-1,\   k \le |\al|-1, \   l \le j
  }}   \mathcal{I}^{r,k,l} \lt| \pl_t^{m-r} \pl^{|\al|-1-k}\bar\pl^{j-l}\pl \oa\rt| , \notag
\end{align}
which,   together with \ef{commutator2},  \ef{2.26-3}, \ef{2.26-4} and \ef{hard}, points to that
\begin{align*}
&(1+t)^{2m+\ka} \int \sigma^{  \iota + |\alpha| } |\mathcal{R}_{2,i}^{m,\al,\ba} |^2 dy \notag\\
   & \ \les
  |\ba|  \sum_{  j\le |\ba|-1} (
 \sum_{r+k+l\le m+|\al|+j-1}\ea_0^2\mathscr{E}^{r+1,k,l}+  \mathscr{E}^{m,|\al|,j}+  \mathscr{E}^{m,|\al|+1,j} )(t) \notag\\
 & \ +  |\alpha|  \sum_{  j\le |\ba|+1} (
 \sum_{r+k+l\le m+|\al|+j-3}\ea_0^2\mathscr{E}^{r+1,k,l}+  \mathscr{E}^{m,|\al|-1,j} )(t)  .
 \end{align*}
This leads us to that for any $\vea>0$,
\begin{subequations}\label{3.14-4}
\begin{align}
&(1+t)^b   \int  \sigma^{  \iota + |\alpha|}  |\mathcal{R}_{2,i}^{m, \al,\ba}  \pl_t^{m+1}  \pl^\alpha \bar\pl^\beta \oa^i | dy\notag\\
& \les   \vea |\ba|   \mathcal{D}_1^{m,\al,\ba}(t)
 +  \vea^{-1} |\ba|  (1+t)^{b-1-2m-\ka}
( \mathscr{E}^{m,|\al|+1,  |\ba|-1}
\notag\\
&+
  \sum_{r+k+l\le m+|\al|+|\ba|-2}\mathscr{E}^{r+1,k,l} )(t), \ \
  \ \  |\al|=0,
\\
&(1+t)^b   \int  \sigma^{  \iota + |\alpha|}  |\mathcal{R}_{2,i}^{m, \al,\ba}  \pl_t^{m+1}  \pl^\alpha \bar\pl^\beta \oa^i | dy\notag\\
& \les \vea   (1+t)^{2m+\kappa +1} \int
 \sigma^{\iota+ |\alpha|}|\pl_t^{m+1}  \pl^\alpha \bar\pl^\beta \oa |^2 dy
\notag\\
& +\vea^{-1}   (1+t)^{2b-4m-2\ka-1}
\sum_{r+k+l\le m+|\al|+|\ba|-1}\mathscr{E}^{r+1,k,l}(t)\notag\\
& \les   \vea (1+t)^{-1} \mathscr{E}^{m+1,|\al|-1, |\ba|}(t)
\notag\\
&+\vea^{-1}  (1+t)^{2b-4m-2\ka-1}
\sum_{r+k+l\le m+|\al|+|\ba|-1}\mathscr{E}^{r+1,k,l}(t), \ \  |\al|\ge 1. \label{5.18-1}
\end{align}
\end{subequations}

As for the term that includes $\mathcal{R}_{3,i}^{m,\al,\ba}$,
it follows from the Taylor expansion that
\begin{align*}
Q(t) & = (1+t)^{-\la}\ta^{n\ga-n+2}-(1+t)\\
& =(1+t)^{-\la} \lt( (\nu+h)^{n\ga-n+2}-\nu^{n\ga-n+2}\rt)\\
& =(1+t)^{-\la}(n\ga-n+2) (\nu+ \vea  h)^{n\ga-n+1}  h
\end{align*}
for some $\vea \in (0,1)$, which  implies,  using \ef{5.26} and \ef{5.27}, that
\be\label{6.20}
|Q(t)|\les \lt\{
\begin{split}
&\ln (1+t), & \la=0,\\
&(1+t)^{\la}, & \la>0.
\end{split}\rt.
\ee
For the derivatives of $Q(t)$, it follows from \ef{pomt} that
$$ (1+t)^{-\la} (\ta^{n\ga -n +2})_t= 1+\la - (n\ga -n +2) \ta^{n\ga -n +1}\ta_{tt}, $$
which gives,
\begin{align*}
Q_t(t)= &\la (1- (1+t)^{-\la-1}\ta^{n\ga-n+2}  ) - (n\ga -n +2) \ta^{n\ga -n +1}\ta_{tt}  \\
=& -\la (1+t)^{-1} Q(t) - (n\ga -n +2) \ta^{n\ga -n +1}\ta_{tt}.
\end{align*}
Based on this, it is easy to derive from \ef{6.20} and \ef{5.26} that for positive integer $h$,
\be\label{6.20-1}
\lt|\frac{d^h Q(t)}{dt^h}\rt|\les C(h) (1+t)^{\la -h }.
\ee
We apply \ef{6.20-1} and  \ef{5.26} to arrive at
$$|\mathcal{R}_{3,i}^{m,\al,\ba} |
     \les  \sum_{1\le h\le m}  (1+t)^{\la -h}
    |\pl_t^{m-h+1}\pl^\alpha \bar\pl^\beta \oa|,  $$
which yields  that for any $\vea>0$,
\begin{align}
&(1+t)^b   \int  \sigma^{  \iota + |\alpha|}  |\mathcal{R}_{3,i}^{m, \al,\ba}  \pl_t^{m+1}  \pl^\alpha \bar\pl^\beta \oa^i | dy \les \vea  \mathcal{D}_1^{m,\al,\ba}(t)
\notag\\
&  +\vea^{-1} (1+t)^{b-2m-\ka-1+(2\la-2)}  \sum_{1\le h\le m}
     \mathscr{E}^{h,|\al|,|\ba|}(t) . \label{3.17-5}
\end{align}
Thanks to \ef{25.2.17},    \ef{6.7-1c},  \ef{7.12}, \ef{6.20} and \ef{5.26},
the remaining terms within $\mathcal{H}_1^{m,\al,\ba}(t)$ can be bounded as indicated below. It holds that for any $\vea>0$,
\begin{align}\label{3.17-3}
&   (1+t)^{b+\la} (1+ \ln (1+t)) \int \sigma^{  \iota + |\alpha|}
 |  \pl_t^{m+1}\pl^\alpha \bar\pl^\beta \oa  |^2 dy
\notag  \\
&  \les \vea \mathcal{D}_1^{m,\al,\ba} (t) +\vea^{-1}(1+t)^{\la-2}(1+\ln^2(1+t)) \mathcal{E}_1^{m,\al,\ba} (t)
\end{align}
and
\begin{align}\label{3.17-4}
& (1+t)^{b}
 \int  \sigma^{  \iota + |\alpha|+1} |\pl_t\pl \oa| |\pl \pl_t^{m} \pl^\alpha \bar\pl^\beta \oa |^2 dy
 \notag\\
 & +(1+t)^{b-1}
 \int  \sigma^{  \iota + |\alpha|} ( | \pl_t^{m} \pl^\alpha \bar\pl^\beta \oa |^2
 +\sa |\pl \pl_t^{m} \pl^\alpha \bar\pl^\beta \oa |^2   )dy
\notag\\& \les ( \vea (1+t)^{-1}+\vea^{-1}(1+t)^{2b-4m-2\ka-1} ) \mathscr{E}^{m,|\al|,|\ba|}(t) .
\end{align}
Hence, as a consequence of \ef{3.14-2}, \ef{3.14-4} and \ef{3.17-5}-\ef{3.17-4}, we arrive at the conclusion that  for any $\vea>0$,
\begin{subequations}\label{3.18}
\begin{align}
& |\mathcal{H}_1^{m,\al,\ba}(t) | \les
\vea^{-1} |\ba| (1+t)^{b-2m-\ka-1}( \mathscr{E}^{m,|\al|+1,  |\ba|-1}
\notag\\
& + \sum_{r+k+l \le  m+|\al|+|\ba|-2} \mathscr{E}^{r+1,k,l} )(t)  +\mathcal{H}_{1good}^{m,\al,\ba}(t),      \ \  |\al|=0, \label{3.18-a} \\
&  |\mathcal{H}_1^{m,\al,\ba}(t) | \les \vea
  (1+t)^{-1} \mathscr{E}^{m+1,|\al|-1, |\ba|}(t)
  +\mathcal{H}_{1good}^{m,\al,\ba}(t)
 ,  \ \    |\al|\ge 1,
 \label{3.18-b}
\end{align}
\end{subequations}
where
\begin{align*}
&\mathcal{H}_{1good}^{m,\al,\ba}(t)
=\vea  (\mathcal{D}_1^{m,\al,\ba} +
  (1+t)^{-1}  \mathscr{E}^{m ,|\al|,|\ba|} )(t)
   +\vea^{-1}(1+t)^{\la-2}(1
   \notag  \\
 & +\ln^2(1+t))\mathcal{E}_1^{m,\al,\ba} (t)
   +\vea^{-1}(1+t)^{b-1-2m-\ka}
 ((1+t)^{2\la-2}
\notag  \\
 & + (1+t)^{b-2m-\ka} )\sum_{r+k+l \le  m+|\al|+|\ba|-1} \mathscr{E}^{r+1,k,l}(t).
\end{align*}
So far, we can find that relying solely on \ef{5.2} is not enough to obtain the estimates we desire. See \ef{3.14-2} and \ef{3.17-4} for example, in addition to the dissipation on $ \pl_t^{m+1}\pl^\alpha \bar\pl^\beta $ as  shown in $\mathcal{D}_1^{m,\al,\ba}$, the dissipation on    $\pl \pl_t^{m}\pl^\alpha \bar\pl^\beta $ is also needed, which will be given in the next step.

{\it Step 2}.  Let $b$ be the time decay rate that appears in \ef{5.2}, then
we integrate the product of $(1+t)^{b-\la} \sigma^{  \iota + |\alpha|}   \pl_t^{m}  \pl^\alpha \bar\pl^\beta \oa^i$ and \ef{25.2.19}   over $\Omega$  and   take \ef{nab} to  generate
\be\label{5.3}
 \frac{d}{dt} \mathcal{E}_4^{m, \al,\ba}(t) +\sum_{j=2,3}\mathcal{D}_j^{m, \al,\ba}  =  \frac{\ta^{n\ga-n+2}}{  (1+t)^{1+\la}}
 \mathcal{D}_1^{m, \al,\ba}-
 \frac{d}{dt}  \mathcal{F}_3^{m, \al,\ba}(t)+  \mathcal{H }_2^{m, \al,\ba}(t),
\ee
where
\begin{align*}
& \mathcal{E}_4^{m,\al,\ba}(t)  =
\frac{1 }{2} (1+t )^{b+1-\la}
 \int  \sigma^{  \iota + |\alpha|}\lt| \pl_t^{m} \pl^\alpha \bar\pl^\beta \oa\rt|^2 dy ,\\
 & \mathcal{D}_2^{m,\al,\ba}(t)  =
\frac{1 }{2}  (2m+2\ka+\la-1 -b) (1+t )^{b-\la}
 \int  \sigma^{  \iota + |\alpha|}\lt| \pl_t^{m} \pl^\alpha \bar\pl^\beta \oa\rt|^2 dy ,\\
 & \mathcal{D}_3^{m,\al,\ba}(t)  =
(1+t)^{b-\la}   \int  \sigma^{\iota+ |\alpha|+1} {J}^{1-\ga}  (   |\na_\eta \pl_t^m  \pl^\alpha \bar\pl^\beta \oa|^2   + \iota^{-1}    |{\rm div}_\eta  \pl_t^m  \pl^\alpha \bar\pl^\beta \oa|^2      ) dy
,\\
& \mathcal{F}_3^{m,\al,\ba}(t)  =
(1+t)^{b-\la}\ta^{n\ga-n+2}
   \int   \sigma^{  \iota + |\alpha|} ( \pl_t^{m}  \pl^\alpha \bar\pl^\beta \oa^i )  \pl_t^{m+1} \pl^\alpha \bar\pl^\beta \oa_i dy ,\\
&\mathcal{H}_2^{m,\al,\ba}(t) = \{((1+t)^{b-\la}\ta^{n\ga-n+2})_t
  - (1+t)^{b-\la}(\ta^{n\ga-n+1}  ( (1+t)^{-\la}\ta+ 2\ta_t ) - (1+t) \\
  & \ +m  (\ta^{n\ga-n+2})_t  )\}
  \int   \sigma^{  \iota + |\alpha|} ( \pl_t^{m}  \pl^\alpha \bar\pl^\beta \oa^i )  \pl_t^{m+1} \pl^\alpha \bar\pl^\beta \oa_i dy  +(1+t)^{-\la}  (2\mathcal{F}_1^{m,\al,\ba}    \\
   & \ +  \mathcal{F}_2^{m,\al,\ba}  )(t)
   +   (1+t)^{b-\la}  \sum_{j=2,3} \int  \sigma^{  \iota + |\alpha|}  \mathcal{R}_{j,i}^{m, \al,\ba}  \pl_t^{m}  \pl^\alpha \bar\pl^\beta \oa^i dy.
\end{align*}
Firstly, to ensure that $\mathcal{D}_2^{m,\al,\ba}$ remains nonnegative, it is necessary to impose the condition that
\be\label{3.18-1}
  b<  2m+2\ka+\la-1.
 \ee
Indeed, we will first choose $b= 2m+\ka+\la-1$ in \ef{3.17} given later, then $b=2m+\ka$ in \ef{5.2}.

In view of  \ef{3.15-2} and that for any $\vea>0$,
$$
 | \mathcal{F}_3^{m, \al,\ba} (t) |  \les   \varepsilon \mathcal{E}_4^{m, \al,\ba}(t)  +  \varepsilon^{-1}
\mathcal{E}_1^{m, \al,\ba}(t),
$$
we see that there exists a large constant $K $ which  only depend on the parameters of the problem,  but does not depend on the data,
such that
\begin{align*}
& \mathcal{E}_I^{m,\al,\ba} (t) =\sum_{1\le j\le 3} K\mathcal{E}_j^{m,\al,\ba} (t)-K\mathcal{F}_2^{m, \al,\ba}(t)
+\mathcal{E}_4^{m, \al,\ba}(t)
+ \mathcal{F}_3^{m, \al,\ba}(t),\\
& \mathcal{D}_I^{m,\al,\ba} (t)=\lt(K-  \frac{\ta^{n\ga-n+2}}{  (1+t)^{1+\la}} \rt)
 \mathcal{D}_1^{m, \al,\ba} (t) +\sum_{j=2,3}\mathcal{D}_j^{m, \al,\ba} (t)
\end{align*}
satisfy
\begin{subequations}\label{3.16}
\begin{align}
& \sum_{|\al|=i, \ |\ba|=j}\mathcal{E}_I^{m,\al,\ba} (t) \les
(1+t)^{b-2m-\ka}   \mathscr{E}^{m,i,j}(t)
 +(1+t)^{b+1-\la}    \int  \sigma^{  \iota +i}\lt| \pl_t^{m} \pl^i \bar\pl^j \oa\rt|^2 dy
\notag\\
& \ \ \ \  +C (1+t)^{b-2m-\ka}
 \sum_{r+k+l\le m+i+j-2}\mathscr{E}^{r+1,k,l} (t),\label{3.16-a}\\
 & \sum_{|\al|=i, \ |\ba|=j}\mathcal{E}_I^{m,\al,\ba} (t) \gtrsim
(1+t)^{b-2m-\ka}  \mathscr{E}^{m,i,j}(t)
 +(1+t)^{b+1-\la}  \int  \sigma^{  \iota +i}\lt| \pl_t^{m} \pl^i \bar\pl^j \oa\rt|^2 dy
\notag\\
& \ \ \ \  -C (1+t)^{b-2m-\ka}
 \sum_{r+k+l\le m+i+j-2}\mathscr{E}^{r+1,k,l} (t),\label{3.16-b}\\
& \sum_{|\al|=i, \ |\ba|=j}\mathcal{D}_I^{m,\al,\ba} (t)
\gtrsim  (1+t)^{b-2m-\ka-\la}  \mathscr{E}^{m,i,j}(t). \label{3.16-c}
\end{align}
\end{subequations}
We will now work on the following equation, which is a linear combination of \ef{5.2} and \ef{5.3}.
\begin{align}\label{3.17}
\frac{d}{dt} \mathcal{E}_I^{m,\al,\ba} (t) +
\mathcal{D}_I^{m,\al,\ba} (t)   =   K  \frac{d}{dt}  \mathcal{F}_1^{m,\al,\ba}(t) + K
  \mathcal{H}_1^{m,\al,\ba}(t)   + \mathcal{H }_2^{m, \al,\ba}(t) .
\end{align}

The estimate of $\mathcal{H }_2^{m, \al,\ba}$ is quite similar to that of $\mathcal{H }_1^{m, \al,\ba}$; here, we will only point out the differences. It holds that for any $\vea>0$,
 \begin{align*}
  &  (1+t)^{b-\la}    \int  \sigma^{  \iota + |\alpha|}  |  \mathcal{R}_{2,i}^{m, \al,\ba}  \pl_t^{m}  \pl^\alpha \bar\pl^\beta \oa^i | dy\\
 & \les \vea|\ba| \mathcal{D}_2^{m,\al,\ba}(t)+\vea^{-1} |\ba|  (1+t)^{b-\la-2m-\ka}
( \mathscr{E}^{m,|\al|+1,  |\ba|-1}
\notag\\
&+
  \sum_{r+k+l\le m+|\al|+|\ba|-2}\mathscr{E}^{r+1,k,l} )(t), \ \  \ \
|\al|=0,\\
  &  (1+t)^{b-\la}    \int  \sigma^{  \iota + |\alpha|}  |  \mathcal{R}_{2,i}^{m, \al,\ba}  \pl_t^{m}  \pl^\alpha \bar\pl^\beta \oa^i | dy\\
 & \les \vea (1+t)^{2b-4m-2\ka-2\la+1}\sum_{r+k+l\le m+|\al|+|\ba|-1}\mathscr{E}^{r+1,k,l}(t)
\\
& +\vea^{-1}(1+t)^{-1}\mathscr{E}^{m,|\al|-1,|\ba|}(t), \ \  \ \
|\al|\ge 1,
\end{align*}
 \begin{align*}
  &  (1+t)^{b-\la}    \int  \sigma^{  \iota + |\alpha|}  | \mathcal{R}_{3,i}^{m, \al,\ba}  \pl_t^{m}  \pl^\alpha \bar\pl^\beta \oa^i | dy
  \\
  & \les \vea (1+t)^{-1}  \sum_{1\le h\le m}
     \mathscr{E}^{h,|\al|,|\ba|}(t)
     +\vea^{-1}(1+t)^{b-2m-\ka+\la -2} \mathcal{E}_4^{m,\al,\ba} (t),
\end{align*}
and
 \begin{align*}
  & (1+t)^b(1+\ln(1+t))
   \int   \sigma^{  \iota + |\alpha|} | ( \pl_t^{m}  \pl^\alpha \bar\pl^\beta \oa^i )  \pl_t^{m+1} \pl^\alpha \bar\pl^\beta \oa_i | dy \\
   &\les \vea \mathcal{D}_1^{m,\al,\ba} (t) + \vea^{-1} (1+t)^{\la-2}(1+\ln^2(1+t)) \mathcal{E}_4^{m,\al,\ba} (t).
\end{align*}
When this is taken with \ef{3.15-1}, \ef{3.15-2}, \ef{6.20}
and \ef{5.26}, it is enough to show that for any $\vea>0$,
\begin{subequations}\label{3.19}
\begin{align}
& |\mathcal{H}_2^{m,\al,\ba}(t) | \les
 \vea^{-1}|\ba| (1+t)^{b-2m-\ka-\la}\mathscr{E}^{m,|\al|+1,|\ba|-1}(t)
\notag \\& \ \ +
 \vea |\ba| \mathcal{D}_2^{m,\al,\ba} (t) +\mathcal{H}_{2good}^{m,\al,\ba}(t),
  \ \ \ \  |\al|=0, \label{3.19-a} \\
&  |\mathcal{H}_2^{m,\al,\ba}(t) | \les  \vea (1+t)^{2b-4m-2\ka-2\la+1}\sum_{r+k+l\le m+|\al|+|\ba|-1}\mathscr{E}^{r+1,k,l}(t)\notag\\ & \ \ +\vea^{-1}(1+t)^{-1}\mathscr{E}^{m,|\al|-1,|\ba|}(t)
+\mathcal{H}_{2good}^{m,\al,\ba}(t)
 ,  \  \  \ \   |\al|\ge 1,
 \label{3.19-b}
\end{align}
\end{subequations}
where
\begin{align*}
& \mathcal{H}_{2good}^{m,\al,\ba}(t) =
 \vea ( \sum_{j=1,3}\mathcal{D}_j^{m,\al,\ba}+(1+t)^{-1} \mathscr{E}^{m,|\al|,|\ba|} ) (t)
\notag  \\& +\vea^{-1} \{
  (1+t)^{\la-2}(1  +\ln^2(1+t))
 + (1+t)^{b-2m-\ka+\la-2}
\}\mathcal{E}_4^{m,\al,\ba} (t)\notag\\ &
+ (1+t)^{b-2m-\ka+\la-2} ( \sum_{r+k+l\le |\al|+|\ba| }\mathscr{E}^{r+1,k,l}(0)
  + \sum_{r+k+l \le m+|\al|+|\ba|-1 }  \sup_{\tau\in [0,t]} \ea_0^2   \mathscr{E}^{r+1,k,l} (\tau)
  ) \notag\\ &
 +(( \vea^{-1}+1) (1+t)^{b-2m-\ka-\la}
 +\vea(1+t)^{-1} )
  \sum_{r+k+l\le m+|\al|+|\ba|-2}\mathscr{E}^{r+1,k,l}(t).
\end{align*}

{\it Step 3}. Based on the preceding analysis, especially on
\ef{3.15-1},  \ef{3.18} and \ef{3.18-1}-\ef{3.19}, we select
 $b= 2m+\ka+\la-1$ in \ef{3.17} to derive the following estimate that for any $\vea\in (0,1)$,
 \begin{align}
&(1+t)^{\la-1}  \mathscr{E}^{m,i,j}(t)
+(1+t)^{2m+\kappa}  \int  \sigma^{  \iota + i}\lt| \pl_t^{m} \pl^i \bar\pl^j \oa\rt|^2 dy
\notag\\
&+\int_0^t (1+\tau)^{-1} \mathscr{E}^{m,i,j}(\tau)d\tau
\les   \vea\int_0^t (1+\tau)^{-1} \mathscr{E}^{m,i,j}(\tau)d\tau
\notag\\
& + \sum_{z=1,4}\mathfrak{G}_z^{m,i,j}(t)
+\vea^{-1}\sum_{z=2,3}\mathfrak{G}_z^{m,i,j}(t), \label{3.27-1}
\end{align}
where
\begin{align*}
& \mathfrak{G}_1^{m,i,j}(t)= \sum_{r+k+l\le m+i+j-1}\{\mathscr{E}^{r+1,k,l}(0)+ \sup_{\tau\in [0,t]} \ea_0^2   \mathscr{E}^{r+1,k,l} (\tau)  \},\\
& \mathfrak{G}_2^{m,i,j}(t)= \sum_{r+k+l\le m+i+j-2}\{(1+t)^{\la-1}  \mathscr{E}^{r+1,k,l}(t)+ \int_0^t (1+\tau)^{-1} \mathscr{E}^{r+1,k,l}(\tau)d\tau
\},\\
& \mathfrak{G}_3^{m,i,j}(t)= \int_0^t  (1+\tau)^{\la-2}(1+\ln^2(1+\tau))   \{\sum_{r+k+l = m+i+j-1}(1+\tau)^{\la-1} \mathscr{E}^{r+1,k,l}(\tau) \\
& \ \ +  (1+\tau)^{2m+\kappa}  \int  \sigma^{  \iota + i}\lt| \pl_\tau^{m} \pl^i\bar\pl^j \oa\rt|^2 dy\} d\tau,\\
& \mathfrak{G}_4^{m,i,j}(t)=\vea^{-1}j
\int_0^t (1+\tau)^{-1} \mathscr{E}^{m,i+1,j-1} (\tau) d\tau, \ \  i=0,\\
&\mathfrak{G}_4^{m,i,j}(t)=\vea \sum_{r+k+l = m+i+j-1}
\int_0^t (1+\tau)^{-1} \mathscr{E}^{r+1,k,l} (\tau) d\tau, \ \  i\ge 1.
\end{align*}
It should be pointed that the lower-order terms $\mathfrak{G}_2^{m,i,j}$   can be bounded  by the mathematical induction,
$\mathfrak{G}_3^{m,i,j}$
by the  Grownwall inequality, while  $\mathfrak{G}_4^{m,i,j}$ by a suitable linear combination of  $(m,0,j)$ and $(m,1,j-1)$ .

First of all, we choose a suitable small value for $\vea$ in \ef{3.27-1} to obtain that
\begin{align}
&(1+t)^{\la-1}  \mathscr{E}^{m,0,j}(t)
+(1+t)^{2m+\kappa}  \int  \sigma^{  \iota }\lt| \pl_t^{m}   \bar\pl^j \oa\rt|^2 dy
+\int_0^t (1+\tau)^{-1} \mathscr{E}^{m,0,j}(\tau)d\tau
\notag\\
& \les  j \int_0^t (1+\tau)^{-1} \mathscr{E}^{m,1,j-1} (\tau) d\tau
+ \sum_{z=1,2,3}\mathfrak{G}_z^{m,0,j}(t).\label{3.27-2}
\end{align}
Thus, when $j\ge 1$,  there exist  large constants $K^{m,j}\ge 1$ which  only depend on the parameters of the problem,  but does not depend on the data, such that for any $\vea\in (0,1)$,
\begin{align}
& (1+t)^{\la-1} (\mathscr{E}^{m,0,j}+ K^{m,j}\mathscr{E}^{m,1,j-1}  ) (t)
+ (1+t)^{2m+\kappa}  \int  \sigma^{  \iota  } (\lt| \pl_t^{m}   \bar\pl^j \oa\rt|^2
\notag\\&
+ K^{m,j} \sa \lt| \pl_t^{m}  \pl^1 \bar\pl^{j-1} \oa\rt|^2
 ) dy
+\int_0^t (1+\tau)^{-1}(   \mathscr{E}^{m,0,j}+  \frac{1}{2} K^{m,j} \mathscr{E}^{m,1,j-1}) (\tau)d\tau
\notag\\
& \les
 \sum_{z=1,2,3}\mathfrak{G}_z^{m,0,j}(t)
 + K^{m,j}(\sum_{z=1,4}\mathfrak{G}_z^{m,1,j-1}
+\vea^{-1} \sum_{z=2,3}\mathfrak{G}_z^{m,1,j-1})(t).
\label{3.27-3}
\end{align}
Hence,  choosing a suitable small value for $\vea$ in \ef{3.27-1} and \ef{3.27-3}, we achieve  that for $h\ge 2$,
\begin{align}
&(\sum_{\substack{m+i+j=h, \\ m\ge 1, \  i=1}}  K^{m, j+1} +\sum_{\substack{m+i+j=h, \\ m\ge 1, \  i \neq 1}} ) \{(1+t)^{\la-1}  \mathscr{E}^{m,i,j}(t)
\notag\\
&+(1+t)^{2m+\kappa}  \int  \sigma^{  \iota + i}\lt| \pl_t^{m} \pl^i \bar\pl^j \oa\rt|^2 dy
+\frac{1}{4}\int_0^t (1+\tau)^{-1} \mathscr{E}^{m,i,j}(\tau)d\tau\}
\notag\\ &
\les  (\sum_{\substack{m+i+j=h, \\ m\ge 1, \  i=1}}  K^{m, j+1} +\sum_{\substack{m+i+j=h, \\ m\ge 1, \  i \neq 1}} )  \sum_{z=1,2,3} \mathfrak{G}_z^{m,i,j}(t),
\notag
\end{align}
where we have used \ef{3.27-2} to deal with the case of $(m,i,j)=(h,0,0)$ in the above summation,   \ef{3.27-3} to  $(m,0,j)$ and $(m,1,j-1)$, and  \ef{3.27-1} to $(m,i,j)$ for $i\ge 2$.
This, together with the Grownwall inequality and the mathematical induction, implies that for $h\ge 1$,
\begin{align}
& \sum_{m+i+j=h, \ m\ge 1 }  \{(1+t)^{\la-1}  \mathscr{E}^{m,i,j}(t)
+(1+t)^{2m+\kappa}  \int  \sigma^{  \iota + i}\lt| \pl_t^{m} \pl^i \bar\pl^j \oa\rt|^2 dy
\notag\\
&
+ \int_0^t (1+\tau)^{-1} \mathscr{E}^{m,i,j}(\tau)d\tau\}
\les \sum_{m+i+j \le h, \ m\ge 1 }\{\mathscr{E}^{m,i,j}(0)+ \sup_{\tau\in [0,t]} \ea_0^2   \mathscr{E}^{m,i,j} (\tau)  \},
\label{3.27-4}
\end{align}
where the below estimate  has served as the basis for deriving  \ef{3.27-4}.
\begin{align*}
&  (1+t)^{\la-1}  \mathscr{E}^{1,0,0}(t)
+(1+t)^{2+\kappa}  \int  \sigma^{  \iota  }\lt| \pl_t  \oa\rt|^2 dy
\notag\\
&
+ \int_0^t (1+\tau)^{-1} \mathscr{E}^{1,0,0}(\tau)d\tau
\les  \mathscr{E}^{1,0,0}(0)+ \sup_{\tau\in [0,t]} \ea_0^2   \mathscr{E}^{1,0,0} (\tau) .
\end{align*}
In fact, the estimate above can be obtained easily, since
$\mathcal{R}_{1,i}^{m,\al,\ba,k}=0$, $\mathcal{R}_{2,i}^{m,\al,\ba}=0$ in \ef{5.2}
 when $m=1$ and $|\al| =|\ba|=0$.

In a similar but simpler way to deriving \ef{3.27-1}, we can
select $b= 2m+\ka$ in \ef{5.2} to find
\begin{align}
&  \mathscr{E}^{m,i,j}(t)
+\int_0^t (1+\tau)^{2m+\ka +1} \int  \sigma^{  \iota + i}
\lt|  \pl_\tau^{m+1}\pl^i \bar\pl^j \oa \rt|^2 dy d\tau
\notag\\ &
\les \sum_{r+k+l\le m+i+j-1}\{\mathscr{E}^{r+1,k,l}(0)+ \sup_{\tau\in [0,t]} \ea_0^2   \mathscr{E}^{r+1,k,l} (\tau)  +\int_0^t (1+\tau)^{-1} \mathscr{E}^{r+1,k,l} (\tau) d\tau\}\notag\\ &
+
\sum_{r+k+l\le m+i+j-2}   \mathscr{E}^{r+1,k,l}(t)
+\int_0^t (1+\tau)^{\la-2}(1+\ln^2(1+\tau)) \mathscr{E}^{m,i,j}(\tau) d\tau,
\notag
\end{align}
which proves \ef{3.28-2}, using \ef{3.27-4},  the Grownwall inequality, the mathematical induction, the smallness of  $\ea_0$, and
\begin{align*}
&    \mathscr{E}^{1,0,0}(t)
+\int_0^t (1+\tau)^{3+\ka} \int  \sigma^{  \iota }
 |  \pl_\tau^{2} \oa |^2 dy d\tau
\les  \mathscr{E}^{1,0,0}(0) .
\end{align*}
\hfill$\Box$

\section{Proof of Theorem \ref{orig}}\label{sec4}
In this section, we prove Theorem \ref{orig}. Clearly, \ef{4.23-a}  follows from \ef{4.24-3}, \ef{25.2.17}, \ef{5.26}, \ef{5.27}, and
$$x(t,y)-\bar x(t,y)= h(t) y + \ta(t) \oa(t,y)$$
for $(t,y)\in[0,\iy)\times \Oa $.
It follows from \ef{Feb25-1}, \ef{4.24-4}, \ef{4.24-1} and  \ef{4.24-2} that for $(t,y)\in[0,\iy)\times \Oa $,
\begin{align*}
&|\rho(t,x(t,y))-\bar\rho(t, \bar x (t,y))|
=\bar\rho_0(y) |\mathscr{J}^{-1}(t,y)-
\nu^{-n}(t)|
\\
& =\bar\rho_0(y) |\ta^{-n}(t) J^{-1}(t,y)-
\nu^{-n}(t)|\\
& = \bar\rho_0(y) |(\ta^{-n}-\nu^{-n})(t) J^{-1}(t,y)+
\nu^{-n}(t)(J^{-1}(t,y)-1)|
\end{align*}
and
\begin{align*}
&|u(t,x(t,y))-\bar u(t, \bar x (t,y))|=|\pl_t x(t,y)-\ka(1+t)^{\ka-1}y|
\\
& =|\ta_t(t)(\oa(t,y)+y)+\ta(t)\pl_t \oa(t,y)
-\nu_t(t) y|\\
& \le | (\ta_t-\nu_t)(t) y|
+|\ta_t(t)\oa(t,y)|+|\ta(t)\pl_t \oa(t,y)|.
\end{align*}
This, together with
\ef{4.24-3},  \ef{25.2.17}, \ef{7.9}, \ef{6.7-1a}, \ef{5.26} and \ef{5.27},  proves \ef{4.23-b} and \ef{4.23-c}.

\renewcommand{\theequation}{A-\arabic{equation}}
\renewcommand{\thethm}{A-\arabic{thm}}
\setcounter{equation}{0}
\setcounter{thm}{0}
\section*{Appendix}  

{\bf The weighted Sobolev embedding}.
Let $\mathfrak{U}$ be a bounded smooth domain in $\mathbb{R}^n  $, and $d=d(y)=dist(y, \partial  \mathfrak{U})$ be a distance function to the boundary.
 For any positive real number $a$  and nonnegative integer $b$,  we define the  weighted Sobolev space  $H^{a, b}(  \mathfrak{U})$   by
$$ H^{a, b}(\mathfrak{U}) = \lt\{   d^{a/2}f \in L^2(\mathfrak{U}): \ \  \int_\mathfrak{U}    d^a|\pl^k f|^2dy<\infty, \ \  0\le k\le b\rt\}$$
  with the norm
$ \|f\|^2_{H^{a, b}(\mathfrak{U})} = \sum_{k=0}^b \int_\mathfrak{U}    d^a|\pl^k f|^2dy$. Let $H^s( \mathfrak{U} )$ $(s\ge 0)$ be the standard Sobolev space, then  for $b\ge  {a}/{2}$, we have the following embedding of weighted Sobolev spaces (cf. \cite{18'}):
 $ H^{a, b}(\mathfrak{U} )\hookrightarrow H^{b- {a}/{2}}( \mathfrak{U})$
    with the estimate
  \be\label{wsv} \|f\|_{H^{b- {a}/{2}}( \mathfrak{U})} \le \bar C \|f\|_{H^{a, b}(\mathfrak{U} )} \ee
for some constant $\bar C$ only depending on $a$, $ b$ and $\mathfrak{U}$.

{\bf The Hardy inequality}. Let $k>-1$ be a given real number, $\bar\vea$ be a positive constant,  and $f$ be a function satisfying
$
\int_0^{\bar\vea} y^{k+2} \lt(f^2 + |f'|^2\rt) dy < \iy,
$
then it holds that
\be\label{hardy'}
\int_0^{\bar\vea} y^{k } f^2 dy \le \bar C \int_0^{\bar\vea} y^{k+2} \lt(f^2 + |f'|^2\rt) dy
\ee
for a certain constant $\bar C$ only depending on $\bar\vea$ and $k$,
whose proof  can be found  in \cite{18'}. Indeed, \ef{hardy'} is a   general version of the  standard Hardy inequality:
$\int_0^\iy |y^{-1} f|^2 dy \le C \int_0^\iy |f'|^2 dy$.
As a consequence of \ef{hardy'}, we have the following  estimates.
Let $k>-1$ be a given real number,
$\sigma(y) = \bar{A}- \bar{B} |{y}|^2$ with positive constants $\bar A$ and $\bar B$,  $\Omega $ be a ball centered at the origin with the radius $ \sqrt{\bar{A}/\bar{B}} $,  and $f$ be a function satisfying $\int_\Omega \sa^{k+2} (f^2 + |\pl f|^2 )dy <\iy$, then it holds that
\begin{align}\label{hard}
\int_\Omega \sa^k f^2 dy \le \bar C \int_\Omega \sa^{k+2} (f^2 + |\pl f|^2) dy
\end{align}
for some constant $\bar C$ only depending on $k, \bar A, \bar B $.
Indeed, \ef{hard} can be found in Lemma 3.2 of \cite{HZeng}.

{\bf The properties of the solution to the ODE problem \ef{pomt}.}
Let $\ka$ and $\nu(t)$ be defined by \ef{5.26-1} and $0\le \la<1$,  then  the solution to problem \ef{pomt}, $h(t)$, and $\ta(t)=\nu(t)+h(t)$  satisfy that for all $t\ge 0$,
\begin{subequations}\label{5.26}
\begin{align}
& C_1(\ka, \la) \nu(t) \le \ta(t) \le C_2(\ka, \la) \nu(t), \ \  0\le \ta_t(t),
\label{5.26-a}\\
& \lt|\frac{d^m \ta(t)}{d t^m}\rt| \le C_3(m, \ka, \la) (1+t)^{\ka -m}, \ \ m=1,2,3, \cdots, \label{5.26-b}
\end{align}
\end{subequations}
and
\begin{subequations}\label{5.27}
\begin{align}
&  |h(t)|\le C_4(\ka, \la) (1+t)^{\ka+\la -1}\lt(1+ \da_{0\la} \ln(1+t)\rt),  \label{6.10-b}\\
&  |h_t(t)|\le C_5(\ka, \la) (1+t)^{\ka+\la -2}\lt(1+ \da_{0\la} \ln(1+t)\rt),  \label{6.10-d}
\end{align}
\end{subequations}
where $C_i(\beta)$ $(i=1,2,\cdots, 5)$ are  positive constants depending only on quantity $\beta$, and $\da_{0\la}$ is the Kronecker Delta symbol satisfying $\da_{0\la}=1$ if $\la=0$, and $\da_{0\la}=0$ if $\la\neq0$. In fact, the estimates \ef{5.26} and  \ef{5.27}  in the case of $\la=0$ and $n=3$ were proved in  Appendix of \cite{HZ}. However, there are some differences between the cases of $\la=0$ and $0<\la<1$, for example, in the behavior of the upper bound of $|h(t)|$: it is an increasing function for $\la>1-\ka$ but a decreasing function for $\la=0$. So our proof will differ somewhat from that in \cite{HZ}, and we will handle both cases of $\la=0$ and $0<\la<1$ in the same manner. The proof of  \ef{5.26} and \ef{5.27}  consists of the following four steps.

{\it Step 1}. In this step, we prove that $\ta_t>0$ for all $t\ge 0$. For this purpose, we write \ef{pomt} as the following system:
\begin{subequations}\begin{align}
&\ta_{tt} +  (1+t)^{-\la} \ta_t=  \ka \ta^{n-n\ga-1 }, \label{4.25-1}\\
& (\ta, \ta_t)(t=0)=(1, \ka).
\end{align}\end{subequations}
Suppose this is not true, due to $\ta_t(0)=\ka >0$,  then there exists $t_*>0$ such that $\ta_t(t)>0$ on  $[0,t_*)$ and
$\ta_t(t_*)=0$, which means $\ta_{tt}(t_*) \le  0$. Substitute these into \ef{4.25-1} to get $ \ka \ta^{n-n\ga-1 }(t_*) \le 0$.
But $\ta(t_*)> \ta(0)=1 $, because of $\ta_t(t)>0$  on $ [0,t_*)$. It is a contradiction. Hence, it holds that
\be\label{4.24-5}
\ta_t(t)  > 0 \ \ {\rm and}  \ \
\ta(t) > 1 \ \
{\rm for} \ \ t > 0.
\ee
Moreover, it is easy to see that
\be\label{4.25-4}
\nu(t)+\vea h(t)  > 1 \ \  {\rm for} \ \  \vea\in [0,1]  \ \  {\rm and } \ \  t>0.
\ee

{\it Step 2}. In this step, we prove that there exist constants $C(\ka, \la)$ such that for all $t\ge 0$,
\be\label{4.27-2}
|h(t)|
\le \lt\{ \begin{split}
&C(\ka ,\la) , & \la+\ka <1,\\
&C(\ka, \la) \ln (1+t), &  \la+\ka =1, \\
& C(\ka, \la) (1+t)^{\la+\ka -1}, & \la+\ka >1.\end{split} \rt.
\ee
To this end, we write \ef{pomt} as the following system:
\begin{subequations}\begin{align}
& h_t=z, \label{4.25-2}\\
& z_t =   -  (1+t)^{-\la} z+  \ka  \lt( (\nu+h)^{n-n\ga-1 } -\nu^{n-n\ga-1 } \rt)- \nu_{ tt}  , \label{4.25-3}\\
& (h,z)(t=0)=(0,0).
 \end{align}\end{subequations}
Since $\ka\in (0,1)$ and $\nu_{tt}=\ka(\ka-1) (1+t)^{\ka -2}<0$, then $z_t>0$ at $t=0$,  so there exists $t_0>0$ such that
$z>0$ and $h>0$ on $(0,t_0)$.
It follows from \ef{4.25-3} that $z_t+  (1+t)^{-\la} z   \le - \nu_{ tt}$ for $h\ge 0$, which implies that for $h\ge 0$,
\be\label{4.27-1}
(e^{\frac{1}{1-\la}(1+t)^{1-\la}}z)_t
\le  \ka(1-\ka) e^{\frac{1}{1-\la}(1+t)^{1-\la}}(1+t)^{\ka -2}.
\ee

Case 1. We assume that $h \ge 0$ on $(0, \iy)$, then it follows from  \ef{25.2.12}, \ef{4.25-2} and \ef{4.27-1}, that  for $t>0$,
$ z(t)  \le C(\ka ,\la) (1+t)^{\la+\ka -2}$
and
\be\label{4.25-7}
0\le h(t)= \int_0^t z(\tau) d\tau
\le \lt\{ \begin{split}
&C(\ka ,\la) , & \la+\ka <1,\\
&C(\ka, \la) \ln (1+t), &  \la+\ka =1, \\
& C(\ka, \la) (1+t)^{\la+\ka -1}, & \la+\ka >1.\end{split} \rt.
\ee
This proves \ef{4.27-2} in this case.

If  Case 1 does not hold,
there exist $t_1<t_2<t_3$ such that  $h(t_1)>0$, $h(t_2)=0$, $h(t_3)<0$, $h\ge 0$ on $(0, t_1)$, $h>0$ on $(t_1, t_2)$, $h<0$ on $(t_2, t_3)$; $z(t_1)=0$
($h$ reaches its positive local maximum at $t_1$),  $z < 0$ on $(t_1, t_3)$ ($h$ decreases from $t_1$ to $t_3$).
Due to the Taylor expansion,   \ef{4.25-3} can be written as
 \begin{align*}
z_t = & -  { (1+t)^{-\la} } z -    \ka  (n\ga-n+1)(1+t)^{-1-\la }  h
 +   2^{-1} (1+\la) (n\ga \notag\\
 &-n+1)
  (\nu +\vea  h)^{n-n\ga-3} h^2 +\ka  (1-\ka)  (1+t)^{\ka-2 }
 \end{align*}
 for some $\vea \in (0,1)$,
which implies, together with \ef{4.25-2} and \ef{4.25-4},  that for $z\le 0$,
\be\label{4.25-5}
\lt(   \ka  (n\ga-n+1)  h^2
+  (1+t)^{\la +1}z^2  \rt)_t\le 0.
\ee

Case 2. We assume that $z\le 0$ on $(t_1,  \iy)$, then it follows from \ef{4.25-5} that for $t>t_1$,
\be\label{4.25-6}
 h^2(t) \le h^2(t_1), \ \
  z^2\le \ka  (n\ga-n+1) h^2(t_1)  (1+t)^{-\la -1},
\ee
where $h(t_1)$ is bounded from above by \ef{4.25-7}. This, along with \ef{4.25-7}, proves \ef{4.27-2} in this case.

If Case 2 does not hold, there exist $t_5>t_4>t_3$ such that
$z<0$ on $(t_1, t_4)$,  $z(t_4)=0$ ($h$ reaches its negative local minimum at $t_4$), $z>0$ on $(t_4, t_5)$; $h(t_5)<0$.

Case 3. We assume that $h\le 0$ on $(t_4, \iy)$, then there could only be two situations: (I) $z>0$ on $(t_4, \iy)$ or  (II) there exist $t_7>t_6>t_5$ such that $z>0$ on $(t_4, t_6)$, $z(t_6)=0$, $z<0$ on $(t_6, t_7)$.  For (I),  we have
$|h(t)|<|h(t_4)|\le h(t_1)$ on $(t_4, \iy)$.  For (II), one has $h(t_7)<h(t_6)\le 0$ which means $|h(t_6)|<|h(t_7)|$,  but it follows from \ef{4.25-6}
that  $|h(t_7)|\le |h(t_6)|$; it is a contradiction so that (II) is ruled out.   This proves  \ef{4.27-2} in this case, due to \ef{4.25-7}.

If Case 3 does not hold, there exist $t_9>t_8 > t_4$ such that
$h(t_8)=0$, $h(t_9)>0$;
$z>0$ on $(t_4,t_9)$.
Clearly, we have $|h(t)|< |h(t_4)|\le h(t_1)$ on $(t_4, t_8)$.
It follows from \ef{4.25-3} and \ef{4.24-5} that $z_t +  (1+t)^{-\la} z\le 2\ka -\ka^2$ on $[t_4, t_8]$, which implies $0<z(t_8)\le C(\ka, \la)$.
We repeat the previous procedure on $[t_8, \infty)$ to obtain \ef{4.27-2} for all $t>0$.

{\it Step 3}. We prove \ef{5.26} in this step. It follows from \ef{4.27-2} that
$|h(t)|\le C(\ka, \la )\nu(t)$ for all $t\ge 0$, which means $\ta(t)\le \nu(t)+ |h(t)| \le C(\ka,\la) \nu(t)$ for all $t\ge 0$.
In view of \ef{4.27-2}, we see that there exists a large time
$T(\ka,\la)$ depending only on $\ka$ and $\la$ such that $|h(t)|\le 2^{-1} \nu(t)$ for $t\ge T(\ka,\la)$, which implies that $\ta(t)\ge \nu(t)-|h(t)|  \ge 2^{-1} \nu(t)$ for $t\ge T(\ka,\la)$. When $t<T(\ka,\la)$, it follows from \ef{4.24-5} that $\ta(t)>1\ge (\nu(T))^{-1} \nu(t)$. This, together with  \ef{4.24-5}, proves \ef{5.26-a}.

As in \cite{HZ}, we can use a mathematical induction to prove \ef{5.26-b} with slight modifications. We will not include the full proof here, but will demonstrate the case of $m=1$ as an example. It follows from \ef{4.25-1} and \ef{5.26-a} that
\begin{align*}
&(e^{\frac{1}{1-\la}(1+t)^{1-\la}}\ta_t )_t
=  \ka e^{\frac{1}{1-\la}(1+t)^{1-\la}}   \ta^{n-n\ga-1 }\\
& \le \ka C_5^{-1}(\ka,\la)  e^{\frac{1}{1-\la}(1+t)^{1-\la}} (1+t)^{\ka (n-n\ga-1)}\\
& =   \ka C_5^{-1}(\ka,\la)  e^{\frac{1}{1-\la}(1+t)^{1-\la}} (1+t)^{\ka -1 -\la},
\end{align*}
which proves \ef{5.26-b} for $m=1$, due to \ef{25.2.12}.

{\it Step 4}. We prove \ef{5.27} in this step.  First, \ef{6.10-b} follows from \ef{4.27-2} in the case of $\la+\ka >1$. For $\la+\ka \le 1$, it follows from \ef{pomt} and the Taylor expansion that
\begin{align}
&    h_t +(1+\la-\ka)(1+t)^{-1} h \notag\\
& = (2\ka)^{-1}(1+\la)(1+\la-\ka)(1+t)^{\la} (\nu+\vea h)^{n-n\ga -3} h^2   - (1+t)^{\la}  \ta_{ tt} \label{6.12}
 \end{align}
for some $\vea \in (0,1)$, which implies, using \ef{5.26}, that
\be\label{6.11-a}
|h(t)|\le C(\ka,\la) (1+t)^{\ka-\la-1} \int_0^t  \lt( (1+\tau )^{\la -2\ka} h^2(\tau) + (1+\tau)^{2\la-1 }   \rt) d\tau.
\ee
When $\la+\ka =1$, we use \ef{4.27-2} and \ef{6.11-a} to get
\begin{align*}
|h(t)|& \le C(\ka,\la) (1+t)^{-2\la} \int_0^t  (1+\tau)^{2\la-1 } \lt( (1+\tau )^{\la -1} \ln^2(1+\tau) +1   \rt) d\tau\\
& \le C(\ka,\la) (1+t)^{-2\la} \int_0^t  (1+\tau)^{2\la-1 }   d\tau \le C(\ka,\la),
\end{align*}
which proves \ef{6.10-b} in the case of $\la+\ka =1$.
When $\la>0$ and $\la+\ka <1$, we use \ef{4.27-2} and \ef{6.11-a} to obtain
\bee
|h(t)|\le \lt\{
\begin{split}
& C(\ka,\la) (1+t)^{\ka+\la-1}  , &\la+2\ka \ge 1,\\
& C(\ka,\la) (1+t)^{-\ka }, &\la+2\ka < 1.
\end{split}\rt.
\eee
This proves  \ef{6.10-b} in the case of $\la+2\ka \ge 1$.
When $\la>0$ and $\la+2\ka <1$, we substitute this into \ef{6.11-a} to give
\bee
|h(t)|\le \lt\{
\begin{split}
& C(\ka,\la) (1+t)^{\ka+\la-1}  , &\la+4\ka \ge 1,\\
& C(\ka,\la) (1+t)^{-3\ka }, &\la+4\ka < 1.
\end{split}\rt.
\eee
This proves  \ef{6.10-b} in the case of $\la+4\ka \ge 1$.
We repeat this procedure to prove \ef{6.10-b} for all $0<\la<1$. Clearly, \ef{6.10-b}for $\la=0$ can be derived in the same way as that for the case of $\la>0$ and $\la+\ka <1$.
Finally,  \ef{6.10-d} follow from  \ef{6.12},
\ef{5.26},
  and \ef{6.10-b}.

\section*{Acknowledgements}
This research was supported in part by NSFC  Grants  12171267 and 11822107, and New Cornerstone Investigator Program 100001127.

\bibliographystyle{plain}

\noindent Huihui Zeng\\
Department of Mathematical Sciences\\
Tsinghua University\\
Beijing, 100084, China;\\
Email: hhzeng@mail.tsinghua.edu.cn

\end{document}